\documentclass[11pt]{article}
\usepackage{amsmath,amsthm,amsfonts,amssymb,url,graphicx,xcolor,authblk,longtable}
\usepackage{tikz}
\usepackage[title]{appendix}
\usetikzlibrary{calc}
\usetikzlibrary{arrows,decorations.markings}
\usepackage{hyperref}
\hypersetup{
    colorlinks=true,
    linkcolor=black,
    filecolor=black,
    urlcolor=blue,
    citecolor=black,
}

\definecolor{red}{RGB}{220,20,60}

\newcommand{\DTS}{\ensuremath{\mathsf{DTS}}} 
\newcommand{\TTS}{\ensuremath{\mathsf{TTS}}}

\newcommand{\Aut}{\ensuremath{\textsc{Aut}}} 

\makeatletter
\g@addto@macro\bfseries{\boldmath}
\makeatother
\newcommand{\D}[3]{\ensuremath{\mathcal{D}_{#1}#2#3}}
\newcommand{\T}[3]{\ensuremath{\mathcal{T}_{#1}#2#3}}

\begin{document}
\title{Good sequencings for small directed triple systems}
\author[1]{Donald L.\ Kreher}
\author[2]{Douglas R.\ Stinson%
\thanks{D.R.\ Stinson's research is supported by  NSERC discovery grant RGPIN-03882.}}
\author[3]{Shannon Veitch}
\affil[1]{Department of Mathematical Sciences, 
Michigan Technological University 
Houghton, MI 49931,  
U.S.A.}
\affil[2]{David R.\ Cheriton School of Computer Science, University of Waterloo,
Waterloo, Ontario, N2L 3G1, Canada}
\affil[3]{Department of Combinatorics and Optimization, University of Waterloo,
Waterloo, Ontario, N2L 3G1, Canada}


\maketitle

\begin{abstract}
A  \emph{directed triple system of order $v$} (or, \DTS$(v)$) is a decomposition of the complete directed graph $\vec{K_v}$
into transitive triples. An \emph{$\ell$-good sequencing}  of a \DTS$(v)$  
is a permutation of the points of the design, 
say $[x_1 \; \cdots \; x_v]$, such that, for every 
triple $(x,y,z)$ in the design,
it is \emph{not} the case that 
$x = x_i$, $y = x_j$ and $z = x_k$ with $i < j < k$ 
and $k-i+1 \leq \ell$.

In this report we provide a maximum $\ell$-good sequencing for each 
\DTS$(v)$, $v \leq 7$.
\end{abstract}

\section{Introduction}
\label{intro.sec}
It is in~\cite{KSV} that the 
\emph{$\ell$-good sequencings} of directed triple systems are first 
introduced and this paper is where we direct the interested reader to find definitions, fundamental results and constructions. In particular it is well known that a \DTS$(v)$ exists if and only if 
$v \equiv 0 \text{ or } 1 \pmod{3}$.
\section{Methodology}

If each triple $(x,y,z)$ in a directed triple system is replaced with the 
3-element subset $\{x,y,z\}$, a twofold triple system (\TTS) is obtained. Thus all \DTS$(v)$
can be obtained by directing each of the triples of a \TTS$(v)$ in all possible ways.
\"{O}stag\r{a}rd and Pottone in~\cite{OP} show that two directed triple systems with nonisomorphic
underlying twofold triple systems cannot be isomorphic. They further show that the automorphism
group of a directed triple system is identical to the automorphism group of its 
underlying twofold triple system. Thus after directing the triples of a \TTS$(v)$
in all possible ways one need only keep the directed triple systems obtained that 
are minimal in their orbit under the action of the automorphism group. Once a \DTS$(v)$ is obtained
we examined all permutations of its vertices in lexicographical order. A permutation
$g$ is an $\ell$-good sequencing if 
there does not exist a triple $(a,b,c)$
such that
\[
g^{-1}(a) < g^{-1}(b) < g^{-1}(c) 
\text{ and }
g^{-1}(c) - g^{-1}(a) \leq \ell -1.
\]
We summarize the results found.
\begin{itemize}
\item There is a unique \TTS$(4)$ and three nonisomorphic $\DTS(4)$. They all  have a 4-good sequencing.
\item There is a unique \TTS$(6)$ and 32 
nonisomorphic $\DTS(6)$. They all  have a 6-good sequencing.
\item There are four nonisomorphic \TTS$(7)$ and 2,368 
nonisomorphic \DTS$(7)$.
Surprisingly  four of the nonisomorphic \DTS$(7)$ do not have a 7-good sequencing. They do have a 6-good sequencing. The remaining
2,364 of the nonisomorphic \DTS$(7)$ have a 6-good sequencing.
\end{itemize}

\section{\boldmath $v=4$}

\subsection*{\boldmath The unique $2$-$(4,3,2)$ design.}

\begin{description}
\item[Design \T{4}{}{}:]
$\{1,2,3\}$,
$\{0,2,3\}$,
$\{0,1,3\}$,
$\{0,1,2\}$.

$\Aut(\T{4}{}{}) = 
\big\langle 
(1,2,3), (0,1)
\big\rangle$, $|\Aut(\T{4}{}{})|=24.$

\end{description}
\subsection{\boldmath The 3 nonisomorphic \DTS$(4)$ with underlying \TTS$(4)$ that is \T{4}{}{} and the good sequencings  found.}
\begin{description}
\item[\D{4}{1}{}:]	$\begin{array}{@{}*{4}{l}@{}}
 (0,3,2)& (1,2,3)&
	 (2,1,0)& (3,0,1)\end{array}$

	Lexicographic least 4-good sequencing 0213

	Number of 4-good sequencings found: 8
\end{description}\begin{description}
\item[\D{4}{2}{}:]$\begin{array}{@{}*{4}{l}@{}}
 (0,3,2)& (1,2,3)&
	 (2,0,1)& (3,1,0)\end{array}$

	Lexicographic least 4-good sequencing found: 0213

	Number of 4-good sequencings found: 8
\end{description}\begin{description}
\item[\D{4}{3}{}:]	$\begin{array}{@{}*{4}{l}@{}}
 (0,3,2)& (1,2,0)&
	 (2,1,3)& (3,0,1)\end{array}$

	Lexicographic least 4-good sequence found: 0123

	Number of 4-good sequencings found: 8
\end{description}

\section{\boldmath $v=6$}

\subsection*{\boldmath The unique $2$-$(6,3,2)$ design.}

\begin{description}
\item[Design \T{6}{}{}:]
$\{2,3,5\}$
$\{2,3,4\}$
$\{1,4,5\}$
$\{1,3,4\}$
$\{1,2,5\}$
$\{0,4,5\}$
$\{0,3,5\}$
$\{0,2,4\}$
$\{0,1,3\}$
$\{0,1,2\}$

$\Aut(\T{6}{}{}) = 
\big\langle 
(1,2,4,5,3),
(0,1)(2,3)
\big\rangle$, $|\Aut(\T{6}{}{})|=24.$

\end{description}

\subsection{\boldmath The 32 nonisomorphic \DTS$(6)$ with underlying \TTS$(6)$ that is \T{6}{}{} and the good sequencings found.}

\begin{description}
\item[\D{6}{1}{}:]	$\begin{array}{@{}*{6}{l}@{}}
 (1,0,3)& (1,4,5)& (2,0,1)& (2,5,3)& (3,0,5)\\
	 (3,2,4)& (4,0,2)& (4,3,1)& (5,0,4)& (5,1,2)\end{array}$

	Lexicographic least 6-good sequencing 012354

	Number of 6-good sequencings found: 70
\end{description}\begin{description}
\item[\D{6}{2}{}:]	$\begin{array}{@{}*{6}{l}@{}}
 (1,0,3)& (1,2,5)& (2,0,1)& (2,4,3)& (3,0,5)\\
	 (3,1,4)& (4,0,2)& (4,5,1)& (5,0,4)& (5,3,2)\end{array}$

	Lexicographic least 6-good sequencing 013452

	Number of 6-good sequencings found: 70
\end{description}\begin{description}
\item[\D{6}{3}{}:]	$\begin{array}{@{}*{6}{l}@{}}
 (0,5,4)& (1,4,5)& (2,1,0)& (2,3,5)& (3,0,1)\\
	 (3,2,4)& (4,0,2)& (4,1,3)& (5,0,3)& (5,1,2)\end{array}$

	Lexicographic least 6-good sequencing 024531

	Number of 6-good sequencings: 37
\end{description}\begin{description}
\item[\D{6}{4}{}:]	$\begin{array}{@{}*{6}{l}@{}}
 (0,5,4)& (1,4,5)& (2,1,0)& (2,3,4)& (3,0,1)\\
	 (3,2,5)& (4,0,2)& (4,1,3)& (5,0,3)& (5,1,2)\end{array}$

	Lexicographic least 6-good sequence found: 024315

	Number of 6-good sequencings found: 42
\end{description}\begin{description}
\item[\D{6}{5}{}:]	$\begin{array}{@{}*{6}{l}@{}}
 (0,5,4)& (1,4,3)& (2,1,0)& (2,3,5)& (3,0,1)\\
	 (3,2,4)& (4,0,2)& (4,1,5)& (5,0,3)& (5,1,2)\end{array}$

	Lexicographic least 6-good sequence found: 013425

	Number of 6-good sequencings found: 43
\end{description}\begin{description}
\item[\D{6}{6}{}:]	$\begin{array}{@{}*{6}{l}@{}}
 (0,5,4)& (1,4,3)& (2,1,0)& (2,3,4)& (3,0,1)\\
	 (3,2,5)& (4,0,2)& (4,1,5)& (5,0,3)& (5,1,2)\end{array}$

	Lexicographic least 6-good sequence found: 013452

	Number of 6-good sequencings found: 42
\end{description}\begin{description}
\item[\D{6}{7}{}:]	$\begin{array}{@{}*{6}{l}@{}}
 (0,5,4)& (1,4,3)& (2,0,1)& (2,3,5)& (3,1,0)\\
	 (3,2,4)& (4,0,2)& (4,1,5)& (5,0,3)& (5,1,2)\end{array}$

	Lexicographic least 6-good sequence found: 013425

	Number of 6-good sequencings found: 48
\end{description}\begin{description}
\item[\D{6}{8}{}:]	$\begin{array}{@{}*{6}{l}@{}}
 (0,5,4)& (1,4,3)& (2,0,1)& (2,3,4)& (3,1,0)\\
	 (3,2,5)& (4,0,2)& (4,1,5)& (5,0,3)& (5,1,2)\end{array}$

	Lexicographic least 6-good sequence found: 013452

	Number of 6-good sequencings found: 40
\end{description}\begin{description}
\item[\D{6}{9}{}:]	$\begin{array}{@{}*{6}{l}@{}}
 (0,5,4)& (1,4,3)& (1,5,2)& (2,1,0)& (2,3,5)\\
	 (3,0,1)& (3,2,4)& (4,0,2)& (4,5,1)& (5,0,3)\end{array}$

	Lexicographic least 6-good sequence found: 013425

	Number of 6-good sequencings found: 43
\end{description}\begin{description}
\item[\D{6}{10}{}:]	$\begin{array}{@{}*{6}{l}@{}}
 (0,5,4)& (1,4,3)& (1,5,2)& (2,1,0)& (2,3,4)\\
	 (3,0,1)& (3,2,5)& (4,0,2)& (4,5,1)& (5,0,3)\end{array}$

	Lexicographic least 6-good sequence found: 024135

	Number of 6-good sequencings found: 42
\end{description}\begin{description}
\item[\D{6}{11}{}:]	$\begin{array}{@{}*{6}{l}@{}}
 (0,5,4)& (1,4,3)& (1,5,2)& (2,0,1)& (2,3,5)\\
	 (3,1,0)& (3,2,4)& (4,0,2)& (4,5,1)& (5,0,3)\end{array}$

	Lexicographic least 6-good sequence found: 013425

	Number of 6-good sequencings found: 48
\end{description}\begin{description}
\item[\D{6}{12}{}:]	$\begin{array}{@{}*{6}{l}@{}}
 (0,5,4)& (1,4,3)& (1,5,2)& (2,0,1)& (2,3,4)\\
	 (3,1,0)& (3,2,5)& (4,0,2)& (4,5,1)& (5,0,3)\end{array}$

	Lexicographic least 6-good sequence found: 024135

	Number of 6-good sequencings found: 45
\end{description}\begin{description}
\item[\D{6}{13}{}:]	$\begin{array}{@{}*{6}{l}@{}}
 (0,5,4)& (1,3,4)& (2,1,0)& (2,4,3)& (3,0,1)\\
	 (3,2,5)& (4,0,2)& (4,1,5)& (5,0,3)& (5,1,2)\end{array}$

	Lexicographic least 6-good sequence found: 014235

	Number of 6-good sequencings found: 43
\end{description}\begin{description}
\item[\D{6}{14}{}:]	$\begin{array}{@{}*{6}{l}@{}}
 (0,5,4)& (1,3,4)& (2,0,1)& (2,4,3)& (3,1,0)\\
	 (3,2,5)& (4,0,2)& (4,1,5)& (5,0,3)& (5,1,2)\end{array}$

	Lexicographic least 6-good sequence found: 014235

	Number of 6-good sequencings found: 42
\end{description}\begin{description}
\item[\D{6}{15}{}:]	$\begin{array}{@{}*{6}{l}@{}}
 (0,5,4)& (1,3,4)& (1,5,2)& (2,1,0)& (2,4,3)\\
	 (3,0,1)& (3,2,5)& (4,0,2)& (4,5,1)& (5,0,3)\end{array}$

	Lexicographic least 6-good sequence found: 014235

	Number of 6-good sequencings found: 43
\end{description}\begin{description}
\item[\D{6}{16}{}:]	$\begin{array}{@{}*{6}{l}@{}}
 (0,5,4)& (1,3,4)& (1,5,2)& (2,0,1)& (2,4,3)\\
	 (3,1,0)& (3,2,5)& (4,0,2)& (4,5,1)& (5,0,3)\end{array}$

	Lexicographic least 6-good sequence found: 014235

	Number of 6-good sequencings found: 43
\end{description}\begin{description}
\item[\D{6}{17}{}:]	$\begin{array}{@{}*{6}{l}@{}}
 (0,5,4)& (1,3,0)& (2,0,1)& (2,4,3)& (3,1,4)\\
	 (3,2,5)& (4,0,2)& (4,1,5)& (5,0,3)& (5,1,2)\end{array}$

	Lexicographic least 6-good sequence found: 012345

	Number of 6-good sequencings found: 49
\end{description}\begin{description}
\item[\D{6}{18}{}:]	$\begin{array}{@{}*{6}{l}@{}}
 (0,5,4)& (1,3,0)& (1,5,2)& (2,0,1)& (2,4,3)\\
	 (3,1,4)& (3,2,5)& (4,0,2)& (4,5,1)& (5,0,3)\end{array}$

	Lexicographic least 6-good sequence found: 012345

	Number of 6-good sequencings found: 43
\end{description}\begin{description}
\item[\D{6}{19}{}:]	$\begin{array}{@{}*{6}{l}@{}}
 (0,5,4)& (1,3,0)& (1,4,5)& (2,0,1)& (2,4,3)\\
	 (3,2,5)& (3,4,1)& (4,0,2)& (5,0,3)& (5,1,2)\end{array}$

	Lexicographic least 6-good sequence found: 041235

	Number of 6-good sequencings found: 49
\end{description}\begin{description}
\item[\D{6}{20}{}:]	$\begin{array}{@{}*{6}{l}@{}}
 (0,5,4)& (1,3,0)& (1,4,5)& (2,0,1)& (2,3,4)\\
	 (3,2,5)& (4,0,2)& (4,3,1)& (5,0,3)& (5,1,2)\end{array}$

	Lexicographic least 6-good sequence found: 024135

	Number of 6-good sequencings found: 43
\end{description}\begin{description}
\item[\D{6}{21}{}:]	$\begin{array}{@{}*{6}{l}@{}}
 (0,5,4)& (1,2,5)& (1,4,3)& (2,1,0)& (2,3,4)\\
	 (3,0,1)& (3,5,2)& (4,0,2)& (4,5,1)& (5,0,3)\end{array}$

	Lexicographic least 6-good sequence found: 024135

	Number of 6-good sequencings found: 42
\end{description}\begin{description}
\item[\D{6}{22}{}:]	$\begin{array}{@{}*{6}{l}@{}}
 (0,5,4)& (1,2,5)& (1,4,3)& (2,0,1)& (2,3,4)\\
	 (3,1,0)& (3,5,2)& (4,0,2)& (4,5,1)& (5,0,3)\end{array}$

	Lexicographic least 6-good sequence found: 024135

	Number of 6-good sequencings found: 48
\end{description}\begin{description}
\item[\D{6}{23}{}:]	$\begin{array}{@{}*{6}{l}@{}}
 (0,5,4)& (1,2,5)& (1,3,4)& (2,1,0)& (2,4,3)\\
	 (3,0,1)& (3,5,2)& (4,0,2)& (4,5,1)& (5,0,3)\end{array}$

	Lexicographic least 6-good sequence found: 014523

	Number of 6-good sequencings found: 37
\end{description}\begin{description}
\item[\D{6}{24}{}:]	$\begin{array}{@{}*{6}{l}@{}}
 (0,5,4)& (1,2,5)& (1,3,4)& (2,0,1)& (2,4,3)\\
	 (3,1,0)& (3,5,2)& (4,0,2)& (4,5,1)& (5,0,3)\end{array}$

	Lexicographic least 6-good sequence found: 014523

	Number of 6-good sequencings found: 40
\end{description}\begin{description}
\item[\D{6}{25}{}:]	$\begin{array}{@{}*{6}{l}@{}}
 (0,5,4)& (1,2,5)& (1,3,0)& (2,0,1)& (2,4,3)\\
	 (3,1,4)& (3,5,2)& (4,0,2)& (4,5,1)& (5,0,3)\end{array}$

	Lexicographic least 6-good sequence found: 014523

	Number of 6-good sequencings found: 42
\end{description}\begin{description}
\item[\D{6}{26}{}:]	$\begin{array}{@{}*{6}{l}@{}}
 (0,5,4)& (1,2,0)& (1,4,5)& (2,3,4)& (3,0,1)\\
	 (3,2,5)& (4,0,2)& (4,1,3)& (5,0,3)& (5,2,1)\end{array}$

	Lexicographic least 6-good sequence found: 024315

	Number of 6-good sequencings found: 43
\end{description}\begin{description}
\item[\D{6}{27}{}:]	$\begin{array}{@{}*{6}{l}@{}}
 (0,5,4)& (1,2,0)& (1,4,3)& (2,3,4)& (3,0,1)\\
	 (3,2,5)& (4,0,2)& (4,1,5)& (5,0,3)& (5,2,1)\end{array}$

	Lexicographic least 6-good sequence found: 013452

	Number of 6-good sequencings found: 45
\end{description}\begin{description}
\item[\D{6}{28}{}:]	$\begin{array}{@{}*{6}{l}@{}}
 (0,5,4)& (1,2,0)& (1,3,4)& (2,4,3)& (3,0,1)\\
	 (3,2,5)& (4,0,2)& (4,1,5)& (5,0,3)& (5,2,1)\end{array}$

	Lexicographic least 6-good sequence found: 014235

	Number of 6-good sequencings found: 48
\end{description}\begin{description}
\item[\D{6}{29}{}:]	$\begin{array}{@{}*{6}{l}@{}}
 (0,5,4)& (1,0,3)& (2,0,1)& (2,4,3)& (3,1,4)\\
	 (3,2,5)& (4,0,2)& (4,1,5)& (5,1,2)& (5,3,0)\end{array}$

	Lexicographic least 6-good sequence found: 012345

	Number of 6-good sequencings found: 54
\end{description}\begin{description}
\item[\D{6}{30}{}:]	$\begin{array}{@{}*{6}{l}@{}}
 (0,5,4)& (1,0,3)& (1,5,2)& (2,0,1)& (2,4,3)\\
	 (3,1,4)& (3,2,5)& (4,0,2)& (4,5,1)& (5,3,0)\end{array}$

	Lexicographic least 6-good sequence found: 012345

	Number of 6-good sequencings found: 53
\end{description}\begin{description}
\item[\D{6}{31}{}:]	$\begin{array}{@{}*{6}{l}@{}}
 (0,5,4)& (1,0,3)& (1,2,5)& (2,0,1)& (2,4,3)\\
	 (3,1,4)& (3,5,0)& (4,0,2)& (4,5,1)& (5,3,2)\end{array}$

	Lexicographic least 6-good sequence found: 013452

	Number of 6-good sequencings found: 54
\end{description}\begin{description}
\item[\D{6}{32}{}:]	$\begin{array}{@{}*{6}{l}@{}}
 (0,5,4)& (1,0,2)& (1,4,3)& (2,3,5)& (2,4,0)\\
	 (3,0,1)& (3,4,2)& (4,1,5)& (5,0,3)& (5,2,1)\end{array}$

	Lexicographic least 6-good sequence found: 013245

	Number of 6-good sequencings found: 53
\end{description}

\section{\boldmath $v=7$}

\subsection*{\boldmath The 4 non-isomorphic 2-(7,3,2) designs .}

\begin{description}
\item[Design \T{7}{1}{}:]
$\{{2,4,5}\}$
$\{{2,4,5}\}$
$\{{2,3,6}\}$
$\{{2,3,6}\}$
$\{{1,4,6}\}$
$\{{1,4,6}\}$
$\{{1,3,5}\}$
$\{{1,3,5}\}$
$\{{0,5,6}\}$
$\{{0,5,6}\}$
$\{{0,3,4}\}$
$\{{0,3,4}\}$
$\{{0,1,2}\}$
$\{{0,1,2}\}$

$\Aut(\T{7}{1}{}) =
\big\langle
(0,1,3,2,5,6,4),
(3,4)(5,6)
\big\rangle$, $|\Aut(\T{7}{1}{})|=168.$

\item[Design \T{7}{2}{}:]
$\{{2,4,5}\}$
$\{{2,4,5}\}$
$\{{2,3,6}\}$
$\{{2,3,6}\}$
$\{{1,5,6}\}$
$\{{1,4,6}\}$
$\{{1,3,5}\}$
$\{{1,3,4}\}$
$\{{0,5,6}\}$
$\{{0,4,6}\}$
$\{{0,3,5}\}$
$\{{0,3,4}\}$
$\{{0,1,2}\}$
$\{{0,1,2}\}$

$\Aut(\T{7}{2}{}) =
\big\langle
(0,6,5,1,3,4)
,
(3,4)(5,6)
\big\rangle$, $
|\Aut(\T{7}{2}{})|=
48.$

\item[Design \T{7}{3}{}:]
$\{{2,4,6}\}$
$\{{2,4,5}\}$
$\{{2,3,6}\}$
$\{{2,3,5}\}$
$\{{1,5,6}\}$
$\{{1,4,5}\}$
$\{{1,3,6}\}$
$\{{1,3,4}\}$
$\{{0,5,6}\}$
$\{{0,4,6}\}$
$\{{0,3,5}\}$
$\{{0,3,4}\}$
$\{{0,1,2}\}$
$\{{0,1,2}\}$

$\Aut(\T{7}{3}{}) =
\big\langle
(1,2)(4,5)
,
(0,2)(3,4,5,6)
\big\rangle$, $
|\Aut(\T{7}{3}{})|=
24.$

\item[Design \T{7}{4}{}:]
$\{{3,4,5}\}$
$\{{2,5,6}\}$
$\{{2,3,6}\}$
$\{{2,3,4}\}$
$\{{1,4,6}\}$
$\{{1,4,5}\}$
$\{{1,3,6}\}$
$\{{1,2,5}\}$
$\{{0,5,6}\}$
$\{{0,4,6}\}$
$\{{0,3,5}\}$
$\{{0,2,4}\}$
$\{{0,1,3}\}$
$\{{0,1,2}\}$

$\big\langle
(1,2,4,6,5,3)
,
(0,1)(2,3)(4,6)
\big\rangle$, $
|\Aut(\T{7}{4}{})|=42.$

\end{description}

\subsection{\boldmath The 18 nonisomorphic \DTS$(7)$ with underlying \TTS$(7)$ 
that is \T{7}{1}{} and the 
good sequencings found.}

\begin{description}
\item[\D{7}{1}{.1}:]	$
\begin{array}{@{}*{7}{l@{\hspace*{6pt}}}@{}}
(1,0,2)& (2,0,1)& (3,0,4)& (3,1,5)& (3,2,6)& (4,0,3)& (4,1,6)\\
(4,2,5)& (5,0,6)& (5,1,3)& (5,2,4)& (6,0,5)& (6,1,4)& (6,2,3)
\end{array}$

	Lexicographic least 7-good sequencing 0123456

	Number of 7-good sequencings found: 480
\end{description}\begin{description}
\item[\D{7}{1}{.2}:]	$\begin{array}{@{}*{7}{l@{\hspace*{6pt}}}@{}}
 (1,0,2)& (2,0,1)& (2,6,3)& (3,0,4)& (3,1,5)& (3,6,2)& (4,0,3)\\
	 (4,1,6)& (4,2,5)& (5,0,6)& (5,1,3)& (5,2,4)& (6,0,5)& (6,1,4)\end{array}$

	Lexicographic least 7-good sequencing 0123456

	Number of 7-good sequencings found: 360
\end{description}\begin{description}
\item[\D{7}{1}{.3}:]	$\begin{array}{@{}*{7}{l@{\hspace*{6pt}}}@{}}
 (1,0,2)& (2,0,1)& (2,5,4)& (2,6,3)& (3,0,4)& (3,1,5)& (3,6,2)\\
	 (4,0,3)& (4,1,6)& (4,5,2)& (5,0,6)& (5,1,3)& (6,0,5)& (6,1,4)\end{array}$

	Lexicographic least 7-good sequencing 0123456

	Number of 7-good sequencings found: 352
\end{description}\begin{description}
\item[\D{7}{1}{.4}:]	$\begin{array}{@{}*{7}{l@{\hspace*{6pt}}}@{}}
 (1,0,2)& (2,0,1)& (2,4,5)& (2,6,3)& (3,0,4)& (3,1,5)& (3,6,2)\\
	 (4,0,3)& (4,1,6)& (5,0,6)& (5,1,3)& (5,4,2)& (6,0,5)& (6,1,4)\end{array}$

	Lexicographic least 7-good sequencing 0123546

	Number of 7-good sequencings found: 352
\end{description}\begin{description}
\item[\D{7}{1}{.5}:]	$\begin{array}{@{}*{7}{l@{\hspace*{6pt}}}@{}}
 (1,0,2)& (1,6,4)& (2,0,1)& (2,5,4)& (3,0,4)& (3,1,5)& (3,2,6)\\
	 (4,0,3)& (4,5,2)& (4,6,1)& (5,0,6)& (5,1,3)& (6,0,5)& (6,2,3)\end{array}$

	Lexicographic least 7-good sequencing 0123456

	Number of 7-good sequencings found: 336
\end{description}\begin{description}
\item[\D{7}{1}{.6}:]	$\begin{array}{@{}*{7}{l@{\hspace*{6pt}}}@{}}
 (1,0,2)& (1,6,4)& (2,0,1)& (2,4,5)& (3,0,4)& (3,1,5)& (3,2,6)\\
	 (4,0,3)& (4,6,1)& (5,0,6)& (5,1,3)& (5,4,2)& (6,0,5)& (6,2,3)\end{array}$

	Lexicographic least 7-good sequencing 0123546

	Number of 7-good sequencings found: 296
\end{description}\begin{description}
\item[\D{7}{1}{.7}:]	$\begin{array}{@{}*{7}{l@{\hspace*{6pt}}}@{}}
 (0,6,5)& (1,0,2)& (2,0,1)& (3,0,4)& (3,1,5)& (3,2,6)& (4,0,3)\\
	 (4,1,6)& (4,2,5)& (5,1,3)& (5,2,4)& (5,6,0)& (6,1,4)& (6,2,3)\end{array}$

	Lexicographic least 7-good sequencing 0123456

	Number of 7-good sequencings found: 352
\end{description}\begin{description}
\item[\D{7}{1}{.8}:]	$\begin{array}{@{}*{7}{l@{\hspace*{6pt}}}@{}}
 (0,6,5)& (1,0,2)& (2,0,1)& (2,6,3)& (3,0,4)& (3,1,5)& (3,6,2)\\
	 (4,0,3)& (4,1,6)& (4,2,5)& (5,1,3)& (5,2,4)& (5,6,0)& (6,1,4)\end{array}$

	Lexicographic least 7-good sequencing 0123456

	Number of 7-good sequencings found: 292
\end{description}\begin{description}
\item[\D{7}{1}{.9}:]	$\begin{array}{@{}*{7}{l@{\hspace*{6pt}}}@{}}
 (0,6,5)& (1,0,2)& (2,0,1)& (2,5,4)& (3,0,4)& (3,1,5)& (3,2,6)\\
	 (4,0,3)& (4,1,6)& (4,5,2)& (5,1,3)& (5,6,0)& (6,1,4)& (6,2,3)\end{array}$

	Lexicographic least 7-good sequencing 0123456

	Number of 7-good sequencings found: 272
\end{description}\begin{description}
\item[\D{7}{1}{.10}:]	$\begin{array}{@{}*{7}{l@{\hspace*{6pt}}}@{}}
 (0,6,5)& (1,0,2)& (2,0,1)& (2,4,5)& (3,0,4)& (3,1,5)& (3,2,6)\\
	 (4,0,3)& (4,1,6)& (5,1,3)& (5,4,2)& (5,6,0)& (6,1,4)& (6,2,3)\end{array}$

	Lexicographic least 7-good sequencing 0123546

	Number of 7-good sequencings found: 284
\end{description}\begin{description}
\item[\D{7}{1}{.11}:]	$\begin{array}{@{}*{7}{l@{\hspace*{6pt}}}@{}}
 (0,6,5)& (1,0,2)& (2,0,1)& (2,3,6)& (3,0,4)& (3,1,5)& (4,0,3)\\
	 (4,1,6)& (4,2,5)& (5,1,3)& (5,2,4)& (5,6,0)& (6,1,4)& (6,3,2)\end{array}$

	Lexicographic least 7-good sequencing 0124563

	Number of 7-good sequencings found: 272
\end{description}\begin{description}
\item[\D{7}{1}{.12}:]	$\begin{array}{@{}*{7}{l@{\hspace*{6pt}}}@{}}
 (0,6,5)& (1,0,2)& (2,0,1)& (2,3,6)& (2,5,4)& (3,0,4)& (3,1,5)\\
	 (4,0,3)& (4,1,6)& (4,5,2)& (5,1,3)& (5,6,0)& (6,1,4)& (6,3,2)\end{array}$

	Lexicographic least 7-good sequencing 0124563

	Number of 7-good sequencings found: 272
\end{description}\begin{description}
\item[\D{7}{1}{.13}:]	$\begin{array}{@{}*{7}{l@{\hspace*{6pt}}}@{}}
 (0,6,5)& (1,0,2)& (2,0,1)& (2,3,6)& (2,4,5)& (3,0,4)& (3,1,5)\\
	 (4,0,3)& (4,1,6)& (5,1,3)& (5,4,2)& (5,6,0)& (6,1,4)& (6,3,2)\end{array}$

	Lexicographic least 7-good sequencing 0125463

	Number of 7-good sequencings found: 264
\end{description}\begin{description}
\item[\D{7}{1}{.14}:]	$\begin{array}{@{}*{7}{l@{\hspace*{6pt}}}@{}}
 (0,6,5)& (1,0,2)& (1,5,3)& (2,0,1)& (2,5,4)& (3,0,4)& (3,2,6)\\
	 (3,5,1)& (4,0,3)& (4,1,6)& (4,5,2)& (5,6,0)& (6,1,4)& (6,2,3)\end{array}$

	Lexicographic least 7-good sequencing 0123456

	Number of 7-good sequencings found: 328
\end{description}\begin{description}
\item[\D{7}{1}{.15}:]	$\begin{array}{@{}*{7}{l@{\hspace*{6pt}}}@{}}
 (0,6,5)& (1,0,2)& (1,5,3)& (2,0,1)& (2,4,5)& (3,0,4)& (3,2,6)\\
	 (3,5,1)& (4,0,3)& (4,1,6)& (5,4,2)& (5,6,0)& (6,1,4)& (6,2,3)\end{array}$

	Lexicographic least 7-good sequencing 0123546

	Number of 7-good sequencings found: 260
\end{description}\begin{description}
\item[\D{7}{1}{.16}:]	$\begin{array}{@{}*{7}{l@{\hspace*{6pt}}}@{}}
 (0,6,5)& (1,0,2)& (1,5,3)& (2,0,1)& (2,3,6)& (3,0,4)& (3,5,1)\\
	 (4,0,3)& (4,1,6)& (4,2,5)& (5,2,4)& (5,6,0)& (6,1,4)& (6,3,2)\end{array}$

	Lexicographic least 7-good sequencing 0132456

	Number of 7-good sequencings found: 244
\end{description}\begin{description}
\item[\D{7}{1}{.17}:]	$\begin{array}{@{}*{7}{l@{\hspace*{6pt}}}@{}}
 (0,6,5)& (1,0,2)& (1,4,6)& (2,0,1)& (2,3,6)& (3,0,4)& (3,1,5)\\
	 (4,0,3)& (4,2,5)& (5,1,3)& (5,2,4)& (5,6,0)& (6,3,2)& (6,4,1)\end{array}$

	Lexicographic least 7-good sequencing 0125634

	Number of 7-good sequencings found: 312
\end{description}\begin{description}
\item[\D{7}{1}{.18}:]	$\begin{array}{@{}*{7}{l@{\hspace*{6pt}}}@{}}
 (0,4,3)& (0,6,5)& (1,0,2)& (1,5,3)& (2,0,1)& (2,3,6)& (3,4,0)\\
	 (3,5,1)& (4,1,6)& (4,2,5)& (5,2,4)& (5,6,0)& (6,1,4)& (6,3,2)\end{array}$

	Lexicographic least 7-good sequencing 0132456

	Number of 7-good sequencings found: 224
\end{description}

\subsection{\boldmath The 274 nonisomorphic \DTS$(7)$ with underlying \TTS$(7)$ that is \T{7}{2}{} and the good sequencings found.}

\begin{description}
\item[\D{7}{2}{.1}:]	$\begin{array}{@{}*{7}{l@{\hspace*{6pt}}}@{}}
 (1,0,2)& (2,0,1)& (3,0,5)& (3,1,4)& (3,2,6)& (4,0,3)& (4,1,6)\\
	 (4,2,5)& (5,0,6)& (5,1,3)& (5,2,4)& (6,0,4)& (6,1,5)& (6,2,3)\end{array}$

	Lexicographic least 7-good sequencing 0123456

	Number of 7-good sequencings found: 416
\end{description}\begin{description}
\item[\D{7}{2}{.2}:]	$\begin{array}{@{}*{7}{l@{\hspace*{6pt}}}@{}}
 (1,0,2)& (2,0,1)& (2,6,3)& (3,0,5)& (3,1,4)& (3,6,2)& (4,0,3)\\
	 (4,1,6)& (4,2,5)& (5,0,6)& (5,1,3)& (5,2,4)& (6,0,4)& (6,1,5)\end{array}$

	Lexicographic least 7-good sequencing 0123456

	Number of 7-good sequencings found: 304
\end{description}\begin{description}
\item[\D{7}{2}{.3}:]	$\begin{array}{@{}*{7}{l@{\hspace*{6pt}}}@{}}
 (1,0,2)& (2,0,1)& (2,5,4)& (2,6,3)& (3,0,5)& (3,1,4)& (3,6,2)\\
	 (4,0,3)& (4,1,6)& (4,5,2)& (5,0,6)& (5,1,3)& (6,0,4)& (6,1,5)\end{array}$

	Lexicographic least 7-good sequencing 0123456

	Number of 7-good sequencings found: 264
\end{description}\begin{description}
\item[\D{7}{2}{.4}:]	$\begin{array}{@{}*{7}{l@{\hspace*{6pt}}}@{}}
 (1,0,2)& (1,6,5)& (2,0,1)& (3,0,5)& (3,1,4)& (3,2,6)& (4,0,3)\\
	 (4,2,5)& (4,6,1)& (5,0,6)& (5,1,3)& (5,2,4)& (6,0,4)& (6,2,3)\end{array}$

	Lexicographic least 7-good sequencing 0123456

	Number of 7-good sequencings found: 312
\end{description}\begin{description}
\item[\D{7}{2}{.5}:]	$\begin{array}{@{}*{7}{l@{\hspace*{6pt}}}@{}}
 (1,0,2)& (1,6,5)& (2,0,1)& (2,6,3)& (3,0,5)& (3,1,4)& (3,6,2)\\
	 (4,0,3)& (4,2,5)& (4,6,1)& (5,0,6)& (5,1,3)& (5,2,4)& (6,0,4)\end{array}$

	Lexicographic least 7-good sequencing 0123456

	Number of 7-good sequencings found: 304
\end{description}\begin{description}
\item[\D{7}{2}{.6}:]	$\begin{array}{@{}*{7}{l@{\hspace*{6pt}}}@{}}
 (1,0,2)& (1,6,5)& (2,0,1)& (2,5,4)& (3,0,5)& (3,1,4)& (3,2,6)\\
	 (4,0,3)& (4,5,2)& (4,6,1)& (5,0,6)& (5,1,3)& (6,0,4)& (6,2,3)\end{array}$

	Lexicographic least 7-good sequencing 0123456

	Number of 7-good sequencings found: 248
\end{description}\begin{description}
\item[\D{7}{2}{.7}:]	$\begin{array}{@{}*{7}{l@{\hspace*{6pt}}}@{}}
 (1,0,2)& (1,6,5)& (2,0,1)& (2,5,4)& (2,6,3)& (3,0,5)& (3,1,4)\\
	 (3,6,2)& (4,0,3)& (4,5,2)& (4,6,1)& (5,0,6)& (5,1,3)& (6,0,4)\end{array}$

	Lexicographic least 7-good sequencing 0123456

	Number of 7-good sequencings found: 280
\end{description}\begin{description}
\item[\D{7}{2}{.8}:]	$\begin{array}{@{}*{7}{l@{\hspace*{6pt}}}@{}}
 (1,0,2)& (1,6,5)& (2,0,1)& (2,4,5)& (3,0,5)& (3,1,4)& (3,2,6)\\
	 (4,0,3)& (4,6,1)& (5,0,6)& (5,1,3)& (5,4,2)& (6,0,4)& (6,2,3)\end{array}$

	Lexicographic least 7-good sequencing 0123546

	Number of 7-good sequencings found: 248
\end{description}\begin{description}
\item[\D{7}{2}{.9}:]	$\begin{array}{@{}*{7}{l@{\hspace*{6pt}}}@{}}
 (1,0,2)& (1,6,5)& (2,0,1)& (2,4,5)& (2,6,3)& (3,0,5)& (3,1,4)\\
	 (3,6,2)& (4,0,3)& (4,6,1)& (5,0,6)& (5,1,3)& (5,4,2)& (6,0,4)\end{array}$

	Lexicographic least 7-good sequencing 0123546

	Number of 7-good sequencings found: 280
\end{description}\begin{description}
\item[\D{7}{2}{.10}:]	$\begin{array}{@{}*{7}{l@{\hspace*{6pt}}}@{}}
 (1,0,2)& (1,6,5)& (2,0,1)& (2,3,6)& (3,0,5)& (3,1,4)& (4,0,3)\\
	 (4,2,5)& (4,6,1)& (5,0,6)& (5,1,3)& (5,2,4)& (6,0,4)& (6,3,2)\end{array}$

	Lexicographic least 7-good sequencing 0124563

	Number of 7-good sequencings found: 264
\end{description}\begin{description}
\item[\D{7}{2}{.11}:]	$\begin{array}{@{}*{7}{l@{\hspace*{6pt}}}@{}}
 (1,0,2)& (1,6,5)& (2,0,1)& (2,3,6)& (2,5,4)& (3,0,5)& (3,1,4)\\
	 (4,0,3)& (4,5,2)& (4,6,1)& (5,0,6)& (5,1,3)& (6,0,4)& (6,3,2)\end{array}$

	Lexicographic least 7-good sequencing 0124563

	Number of 7-good sequencings found: 232
\end{description}\begin{description}
\item[\D{7}{2}{.12}:]	$\begin{array}{@{}*{7}{l@{\hspace*{6pt}}}@{}}
 (1,0,2)& (1,6,5)& (2,0,1)& (2,3,6)& (2,4,5)& (3,0,5)& (3,1,4)\\
	 (4,0,3)& (4,6,1)& (5,0,6)& (5,1,3)& (5,4,2)& (6,0,4)& (6,3,2)\end{array}$

	Lexicographic least 7-good sequencing 0125463

	Number of 7-good sequencings found: 232
\end{description}\begin{description}
\item[\D{7}{2}{.13}:]	$\begin{array}{@{}*{7}{l@{\hspace*{6pt}}}@{}}
 (1,0,2)& (1,4,6)& (1,5,3)& (2,0,1)& (3,0,5)& (3,2,6)& (3,4,1)\\
	 (4,0,3)& (4,2,5)& (5,0,6)& (5,2,4)& (6,0,4)& (6,2,3)& (6,5,1)\end{array}$

	Lexicographic least 7-good sequencing 0123564

	Number of 7-good sequencings found: 320
\end{description}\begin{description}
\item[\D{7}{2}{.14}:]	$\begin{array}{@{}*{7}{l@{\hspace*{6pt}}}@{}}
 (1,0,2)& (1,4,6)& (1,5,3)& (2,0,1)& (2,6,3)& (3,0,5)& (3,4,1)\\
	 (3,6,2)& (4,0,3)& (4,2,5)& (5,0,6)& (5,2,4)& (6,0,4)& (6,5,1)\end{array}$

	Lexicographic least 7-good sequencing 0123564

	Number of 7-good sequencings found: 264
\end{description}\begin{description}
\item[\D{7}{2}{.15}:]	$\begin{array}{@{}*{7}{l@{\hspace*{6pt}}}@{}}
 (1,0,2)& (1,4,6)& (1,5,3)& (2,0,1)& (2,5,4)& (3,0,5)& (3,2,6)\\
	 (3,4,1)& (4,0,3)& (4,5,2)& (5,0,6)& (6,0,4)& (6,2,3)& (6,5,1)\end{array}$

	Lexicographic least 7-good sequencing 0123645

	Number of 7-good sequencings found: 312
\end{description}\begin{description}
\item[\D{7}{2}{.16}:]	$\begin{array}{@{}*{7}{l@{\hspace*{6pt}}}@{}}
 (1,0,2)& (1,4,6)& (1,5,3)& (2,0,1)& (2,5,4)& (2,6,3)& (3,0,5)\\
	 (3,4,1)& (3,6,2)& (4,0,3)& (4,5,2)& (5,0,6)& (6,0,4)& (6,5,1)\end{array}$

	Lexicographic least 7-good sequencing 0123645

	Number of 7-good sequencings found: 280
\end{description}\begin{description}
\item[\D{7}{2}{.17}:]	$\begin{array}{@{}*{7}{l@{\hspace*{6pt}}}@{}}
 (1,0,2)& (1,4,6)& (1,5,3)& (2,0,1)& (2,4,5)& (2,6,3)& (3,0,5)\\
	 (3,4,1)& (3,6,2)& (4,0,3)& (5,0,6)& (5,4,2)& (6,0,4)& (6,5,1)\end{array}$

	Lexicographic least 7-good sequencing 0123564

	Number of 7-good sequencings found: 280
\end{description}\begin{description}
\item[\D{7}{2}{.18}:]	$\begin{array}{@{}*{7}{l@{\hspace*{6pt}}}@{}}
 (0,6,5)& (1,0,2)& (2,0,1)& (3,0,4)& (3,1,5)& (3,2,6)& (4,1,3)\\
	 (4,2,5)& (4,6,0)& (5,0,3)& (5,1,6)& (5,2,4)& (6,1,4)& (6,2,3)\end{array}$

	Lexicographic least 7-good sequencing 0123456

	Number of 7-good sequencings found: 304
\end{description}\begin{description}
\item[\D{7}{2}{.19}:]	$\begin{array}{@{}*{7}{l@{\hspace*{6pt}}}@{}}
 (0,6,5)& (1,0,2)& (2,0,1)& (2,6,3)& (3,0,4)& (3,1,5)& (3,6,2)\\
	 (4,1,3)& (4,2,5)& (4,6,0)& (5,0,3)& (5,1,6)& (5,2,4)& (6,1,4)\end{array}$

	Lexicographic least 7-good sequencing 0123456

	Number of 7-good sequencings found: 250
\end{description}\begin{description}
\item[\D{7}{2}{.20}:]	$\begin{array}{@{}*{7}{l@{\hspace*{6pt}}}@{}}
 (0,6,5)& (1,0,2)& (2,0,1)& (2,5,4)& (3,0,4)& (3,1,5)& (3,2,6)\\
	 (4,1,3)& (4,5,2)& (4,6,0)& (5,0,3)& (5,1,6)& (6,1,4)& (6,2,3)\end{array}$

	Lexicographic least 7-good sequencing 0123456

	Number of 7-good sequencings found: 222
\end{description}\begin{description}
\item[\D{7}{2}{.21}:]	$\begin{array}{@{}*{7}{l@{\hspace*{6pt}}}@{}}
 (0,6,5)& (1,0,2)& (2,0,1)& (2,5,4)& (2,6,3)& (3,0,4)& (3,1,5)\\
	 (3,6,2)& (4,1,3)& (4,5,2)& (4,6,0)& (5,0,3)& (5,1,6)& (6,1,4)\end{array}$

	Lexicographic least 7-good sequencing 0123456

	Number of 7-good sequencings found: 217
\end{description}\begin{description}
\item[\D{7}{2}{.22}:]	$\begin{array}{@{}*{7}{l@{\hspace*{6pt}}}@{}}
 (0,6,5)& (1,0,2)& (2,0,1)& (2,4,5)& (3,0,4)& (3,1,5)& (3,2,6)\\
	 (4,1,3)& (4,6,0)& (5,0,3)& (5,1,6)& (5,4,2)& (6,1,4)& (6,2,3)\end{array}$

	Lexicographic least 7-good sequencing 0123546

	Number of 7-good sequencings found: 222
\end{description}\begin{description}
\item[\D{7}{2}{.23}:]	$\begin{array}{@{}*{7}{l@{\hspace*{6pt}}}@{}}
 (0,6,5)& (1,0,2)& (2,0,1)& (2,4,5)& (2,6,3)& (3,0,4)& (3,1,5)\\
	 (3,6,2)& (4,1,3)& (4,6,0)& (5,0,3)& (5,1,6)& (5,4,2)& (6,1,4)\end{array}$

	Lexicographic least 7-good sequencing 0123546

	Number of 7-good sequencings found: 217
\end{description}\begin{description}
\item[\D{7}{2}{.24}:]	$\begin{array}{@{}*{7}{l@{\hspace*{6pt}}}@{}}
 (0,6,5)& (1,0,2)& (2,0,1)& (2,3,6)& (3,0,4)& (3,1,5)& (4,1,3)\\
	 (4,2,5)& (4,6,0)& (5,0,3)& (5,1,6)& (5,2,4)& (6,1,4)& (6,3,2)\end{array}$

	Lexicographic least 7-good sequencing 0124563

	Number of 7-good sequencings found: 230
\end{description}\begin{description}
\item[\D{7}{2}{.25}:]	$\begin{array}{@{}*{7}{l@{\hspace*{6pt}}}@{}}
 (0,6,5)& (1,0,2)& (2,0,1)& (2,3,6)& (2,5,4)& (3,0,4)& (3,1,5)\\
	 (4,1,3)& (4,5,2)& (4,6,0)& (5,0,3)& (5,1,6)& (6,1,4)& (6,3,2)\end{array}$

	Lexicographic least 7-good sequencing 0124563

	Number of 7-good sequencings found: 193
\end{description}\begin{description}
\item[\D{7}{2}{.26}:]	$\begin{array}{@{}*{7}{l@{\hspace*{6pt}}}@{}}
 (0,6,5)& (1,0,2)& (2,0,1)& (2,3,6)& (2,4,5)& (3,0,4)& (3,1,5)\\
	 (4,1,3)& (4,6,0)& (5,0,3)& (5,1,6)& (5,4,2)& (6,1,4)& (6,3,2)\end{array}$

	Lexicographic least 7-good sequencing 0125463

	Number of 7-good sequencings found: 193
\end{description}\begin{description}
\item[\D{7}{2}{.27}:]	$\begin{array}{@{}*{7}{l@{\hspace*{6pt}}}@{}}
 (0,6,5)& (1,0,2)& (1,6,4)& (2,0,1)& (3,0,4)& (3,1,5)& (3,2,6)\\
	 (4,1,3)& (4,2,5)& (4,6,0)& (5,0,3)& (5,2,4)& (5,6,1)& (6,2,3)\end{array}$

	Lexicographic least 7-good sequencing 0123456

	Number of 7-good sequencings found: 236
\end{description}\begin{description}
\item[\D{7}{2}{.28}:]	$\begin{array}{@{}*{7}{l@{\hspace*{6pt}}}@{}}
 (0,6,5)& (1,0,2)& (1,6,4)& (2,0,1)& (2,5,4)& (3,0,4)& (3,1,5)\\
	 (3,2,6)& (4,1,3)& (4,5,2)& (4,6,0)& (5,0,3)& (5,6,1)& (6,2,3)\end{array}$

	Lexicographic least 7-good sequencing 0123456

	Number of 7-good sequencings found: 190
\end{description}\begin{description}
\item[\D{7}{2}{.29}:]	$\begin{array}{@{}*{7}{l@{\hspace*{6pt}}}@{}}
 (0,6,5)& (1,0,2)& (1,6,4)& (2,0,1)& (2,4,5)& (3,0,4)& (3,1,5)\\
	 (3,2,6)& (4,1,3)& (4,6,0)& (5,0,3)& (5,4,2)& (5,6,1)& (6,2,3)\end{array}$

	Lexicographic least 7-good sequencing 0123546

	Number of 7-good sequencings found: 190
\end{description}\begin{description}
\item[\D{7}{2}{.30}:]	$\begin{array}{@{}*{7}{l@{\hspace*{6pt}}}@{}}
 (0,6,5)& (1,0,2)& (1,6,4)& (2,0,1)& (2,3,6)& (3,0,4)& (3,1,5)\\
	 (4,1,3)& (4,2,5)& (4,6,0)& (5,0,3)& (5,2,4)& (5,6,1)& (6,3,2)\end{array}$

	Lexicographic least 7-good sequencing 0124563

	Number of 7-good sequencings found: 236
\end{description}\begin{description}
\item[\D{7}{2}{.31}:]	$\begin{array}{@{}*{7}{l@{\hspace*{6pt}}}@{}}
 (0,6,5)& (1,0,2)& (1,6,4)& (2,0,1)& (2,3,6)& (2,5,4)& (3,0,4)\\
	 (3,1,5)& (4,1,3)& (4,5,2)& (4,6,0)& (5,0,3)& (5,6,1)& (6,3,2)\end{array}$

	Lexicographic least 7-good sequencing 0124563

	Number of 7-good sequencings found: 190
\end{description}\begin{description}
\item[\D{7}{2}{.32}:]	$\begin{array}{@{}*{7}{l@{\hspace*{6pt}}}@{}}
 (0,6,5)& (1,0,2)& (1,6,4)& (2,0,1)& (2,3,6)& (2,4,5)& (3,0,4)\\
	 (3,1,5)& (4,1,3)& (4,6,0)& (5,0,3)& (5,4,2)& (5,6,1)& (6,3,2)\end{array}$

	Lexicographic least 7-good sequencing 0125463

	Number of 7-good sequencings found: 190
\end{description}\begin{description}
\item[\D{7}{2}{.33}:]	$\begin{array}{@{}*{7}{l@{\hspace*{6pt}}}@{}}
 (0,6,5)& (1,0,2)& (1,5,6)& (2,0,1)& (3,1,4)& (3,2,6)& (3,5,0)\\
	 (4,0,3)& (4,2,5)& (4,6,1)& (5,1,3)& (5,2,4)& (6,0,4)& (6,2,3)\end{array}$

	Lexicographic least 7-good sequencing 0234516

	Number of 7-good sequencings found: 148
\end{description}\begin{description}
\item[\D{7}{2}{.34}:]	$\begin{array}{@{}*{7}{l@{\hspace*{6pt}}}@{}}
 (0,6,5)& (1,0,2)& (1,5,6)& (2,0,1)& (3,0,4)& (3,2,6)& (3,5,1)\\
	 (4,1,3)& (4,2,5)& (4,6,0)& (5,0,3)& (5,2,4)& (6,1,4)& (6,2,3)\end{array}$

	Lexicographic least 7-good sequencing 0245316

	Number of 7-good sequencings found: 206
\end{description}\begin{description}
\item[\D{7}{2}{.35}:]	$\begin{array}{@{}*{7}{l@{\hspace*{6pt}}}@{}}
 (0,6,5)& (1,0,2)& (1,5,6)& (2,0,1)& (2,6,3)& (3,1,4)& (3,5,0)\\
	 (3,6,2)& (4,0,3)& (4,2,5)& (4,6,1)& (5,1,3)& (5,2,4)& (6,0,4)\end{array}$

	Lexicographic least 7-good sequencing 0234516

	Number of 7-good sequencings found: 166
\end{description}\begin{description}
\item[\D{7}{2}{.36}:]	$\begin{array}{@{}*{7}{l@{\hspace*{6pt}}}@{}}
 (0,6,5)& (1,0,2)& (1,5,6)& (2,0,1)& (2,6,3)& (3,0,4)& (3,5,1)\\
	 (3,6,2)& (4,1,3)& (4,2,5)& (4,6,0)& (5,0,3)& (5,2,4)& (6,1,4)\end{array}$

	Lexicographic least 7-good sequencing 0245316

	Number of 7-good sequencings found: 191
\end{description}\begin{description}
\item[\D{7}{2}{.37}:]	$\begin{array}{@{}*{7}{l@{\hspace*{6pt}}}@{}}
 (0,6,5)& (1,0,2)& (1,5,6)& (2,0,1)& (2,5,4)& (3,1,4)& (3,2,6)\\
	 (3,5,0)& (4,0,3)& (4,5,2)& (4,6,1)& (5,1,3)& (6,0,4)& (6,2,3)\end{array}$

	Lexicographic least 7-good sequencing 0234516

	Number of 7-good sequencings found: 168
\end{description}\begin{description}
\item[\D{7}{2}{.38}:]	$\begin{array}{@{}*{7}{l@{\hspace*{6pt}}}@{}}
 (0,6,5)& (1,0,2)& (1,5,6)& (2,0,1)& (2,5,4)& (3,0,4)& (3,2,6)\\
	 (3,5,1)& (4,1,3)& (4,5,2)& (4,6,0)& (5,0,3)& (6,1,4)& (6,2,3)\end{array}$

	Lexicographic least 7-good sequencing 0245316

	Number of 7-good sequencings found: 210
\end{description}\begin{description}
\item[\D{7}{2}{.39}:]	$\begin{array}{@{}*{7}{l@{\hspace*{6pt}}}@{}}
 (0,6,5)& (1,0,2)& (1,5,6)& (2,0,1)& (2,5,4)& (2,6,3)& (3,1,4)\\
	 (3,5,0)& (3,6,2)& (4,0,3)& (4,5,2)& (4,6,1)& (5,1,3)& (6,0,4)\end{array}$

	Lexicographic least 7-good sequencing 0234516

	Number of 7-good sequencings found: 195
\end{description}\begin{description}
\item[\D{7}{2}{.40}:]	$\begin{array}{@{}*{7}{l@{\hspace*{6pt}}}@{}}
 (0,6,5)& (1,0,2)& (1,5,6)& (2,0,1)& (2,5,4)& (2,6,3)& (3,0,4)\\
	 (3,5,1)& (3,6,2)& (4,1,3)& (4,5,2)& (4,6,0)& (5,0,3)& (6,1,4)\end{array}$

	Lexicographic least 7-good sequencing 0245316

	Number of 7-good sequencings found: 214
\end{description}\begin{description}
\item[\D{7}{2}{.41}:]	$\begin{array}{@{}*{7}{l@{\hspace*{6pt}}}@{}}
 (0,6,5)& (1,0,2)& (1,5,6)& (2,0,1)& (2,4,5)& (3,1,4)& (3,2,6)\\
	 (3,5,0)& (4,0,3)& (4,6,1)& (5,1,3)& (5,4,2)& (6,0,4)& (6,2,3)\end{array}$

	Lexicographic least 7-good sequencing 0235416

	Number of 7-good sequencings found: 122
\end{description}\begin{description}
\item[\D{7}{2}{.42}:]	$\begin{array}{@{}*{7}{l@{\hspace*{6pt}}}@{}}
 (0,6,5)& (1,0,2)& (1,5,6)& (2,0,1)& (2,4,5)& (3,0,4)& (3,2,6)\\
	 (3,5,1)& (4,1,3)& (4,6,0)& (5,0,3)& (5,4,2)& (6,1,4)& (6,2,3)\end{array}$

	Lexicographic least 7-good sequencing 0251346

	Number of 7-good sequencings found: 184
\end{description}\begin{description}
\item[\D{7}{2}{.43}:]	$\begin{array}{@{}*{7}{l@{\hspace*{6pt}}}@{}}
 (0,6,5)& (1,0,2)& (1,5,6)& (2,0,1)& (2,4,5)& (2,6,3)& (3,1,4)\\
	 (3,5,0)& (3,6,2)& (4,0,3)& (4,6,1)& (5,1,3)& (5,4,2)& (6,0,4)\end{array}$

	Lexicographic least 7-good sequencing 0235416

	Number of 7-good sequencings found: 157
\end{description}\begin{description}
\item[\D{7}{2}{.44}:]	$\begin{array}{@{}*{7}{l@{\hspace*{6pt}}}@{}}
 (0,6,5)& (1,0,2)& (1,5,6)& (2,0,1)& (2,4,5)& (2,6,3)& (3,0,4)\\
	 (3,5,1)& (3,6,2)& (4,1,3)& (4,6,0)& (5,0,3)& (5,4,2)& (6,1,4)\end{array}$

	Lexicographic least 7-good sequencing 0251346

	Number of 7-good sequencings found: 181
\end{description}\begin{description}
\item[\D{7}{2}{.45}:]	$\begin{array}{@{}*{7}{l@{\hspace*{6pt}}}@{}}
 (0,6,5)& (1,0,2)& (1,5,6)& (2,0,1)& (2,3,6)& (3,1,4)& (3,5,0)\\
	 (4,0,3)& (4,2,5)& (4,6,1)& (5,1,3)& (5,2,4)& (6,0,4)& (6,3,2)\end{array}$

	Lexicographic least 7-good sequencing 0256341

	Number of 7-good sequencings found: 158
\end{description}\begin{description}
\item[\D{7}{2}{.46}:]	$\begin{array}{@{}*{7}{l@{\hspace*{6pt}}}@{}}
 (0,6,5)& (1,0,2)& (1,5,6)& (2,0,1)& (2,3,6)& (3,0,4)& (3,5,1)\\
	 (4,1,3)& (4,2,5)& (4,6,0)& (5,0,3)& (5,2,4)& (6,1,4)& (6,3,2)\end{array}$

	Lexicographic least 7-good sequencing 0245631

	Number of 7-good sequencings found: 189
\end{description}\begin{description}
\item[\D{7}{2}{.47}:]	$\begin{array}{@{}*{7}{l@{\hspace*{6pt}}}@{}}
 (0,6,5)& (1,0,2)& (1,5,6)& (2,0,1)& (2,3,6)& (2,5,4)& (3,1,4)\\
	 (3,5,0)& (4,0,3)& (4,5,2)& (4,6,1)& (5,1,3)& (6,0,4)& (6,3,2)\end{array}$

	Lexicographic least 7-good sequencing 0324516

	Number of 7-good sequencings found: 179
\end{description}\begin{description}
\item[\D{7}{2}{.48}:]	$\begin{array}{@{}*{7}{l@{\hspace*{6pt}}}@{}}
 (0,6,5)& (1,0,2)& (1,5,6)& (2,0,1)& (2,3,6)& (2,5,4)& (3,0,4)\\
	 (3,5,1)& (4,1,3)& (4,5,2)& (4,6,0)& (5,0,3)& (6,1,4)& (6,3,2)\end{array}$

	Lexicographic least 7-good sequencing 0245631

	Number of 7-good sequencings found: 196
\end{description}\begin{description}
\item[\D{7}{2}{.49}:]	$\begin{array}{@{}*{7}{l@{\hspace*{6pt}}}@{}}
 (0,6,5)& (1,0,2)& (1,5,6)& (2,0,1)& (2,3,6)& (2,4,5)& (3,1,4)\\
	 (3,5,0)& (4,0,3)& (4,6,1)& (5,1,3)& (5,4,2)& (6,0,4)& (6,3,2)\end{array}$

	Lexicographic least 7-good sequencing 0256341

	Number of 7-good sequencings found: 139
\end{description}\begin{description}
\item[\D{7}{2}{.50}:]	$\begin{array}{@{}*{7}{l@{\hspace*{6pt}}}@{}}
 (0,6,5)& (1,0,2)& (1,5,6)& (2,0,1)& (2,3,6)& (2,4,5)& (3,0,4)\\
	 (3,5,1)& (4,1,3)& (4,6,0)& (5,0,3)& (5,4,2)& (6,1,4)& (6,3,2)\end{array}$

	Lexicographic least 7-good sequencing 0251463

	Number of 7-good sequencings found: 149
\end{description}\begin{description}
\item[\D{7}{2}{.51}:]	$\begin{array}{@{}*{7}{l@{\hspace*{6pt}}}@{}}
 (0,6,5)& (1,0,2)& (1,5,3)& (2,0,1)& (3,1,4)& (3,2,6)& (3,5,0)\\
	 (4,0,3)& (4,2,5)& (4,6,1)& (5,1,6)& (5,2,4)& (6,0,4)& (6,2,3)\end{array}$

	Lexicographic least 7-good sequencing 0123456

	Number of 7-good sequencings found: 153
\end{description}\begin{description}
\item[\D{7}{2}{.52}:]	$\begin{array}{@{}*{7}{l@{\hspace*{6pt}}}@{}}
 (0,6,5)& (1,0,2)& (1,5,3)& (2,0,1)& (3,1,4)& (3,2,6)& (3,5,0)\\
	 (4,0,3)& (4,1,6)& (4,2,5)& (5,2,4)& (5,6,1)& (6,0,4)& (6,2,3)\end{array}$

	Lexicographic least 7-good sequencing 0123456

	Number of 7-good sequencings found: 152
\end{description}\begin{description}
\item[\D{7}{2}{.53}:]	$\begin{array}{@{}*{7}{l@{\hspace*{6pt}}}@{}}
 (0,6,5)& (1,0,2)& (1,5,3)& (2,0,1)& (2,6,3)& (3,1,4)& (3,5,0)\\
	 (3,6,2)& (4,0,3)& (4,2,5)& (4,6,1)& (5,1,6)& (5,2,4)& (6,0,4)\end{array}$

	Lexicographic least 7-good sequencing 0123456

	Number of 7-good sequencings found: 155
\end{description}\begin{description}
\item[\D{7}{2}{.54}:]	$\begin{array}{@{}*{7}{l@{\hspace*{6pt}}}@{}}
 (0,6,5)& (1,0,2)& (1,5,3)& (2,0,1)& (2,6,3)& (3,1,4)& (3,5,0)\\
	 (3,6,2)& (4,0,3)& (4,1,6)& (4,2,5)& (5,2,4)& (5,6,1)& (6,0,4)\end{array}$

	Lexicographic least 7-good sequencing 0123456

	Number of 7-good sequencings found: 148
\end{description}\begin{description}
\item[\D{7}{2}{.55}:]	$\begin{array}{@{}*{7}{l@{\hspace*{6pt}}}@{}}
 (0,6,5)& (1,0,2)& (1,5,3)& (2,0,1)& (2,5,4)& (3,1,4)& (3,2,6)\\
	 (3,5,0)& (4,0,3)& (4,5,2)& (4,6,1)& (5,1,6)& (6,0,4)& (6,2,3)\end{array}$

	Lexicographic least 7-good sequencing 0123456

	Number of 7-good sequencings found: 177
\end{description}\begin{description}
\item[\D{7}{2}{.56}:]	$\begin{array}{@{}*{7}{l@{\hspace*{6pt}}}@{}}
 (0,6,5)& (1,0,2)& (1,5,3)& (2,0,1)& (2,5,4)& (3,1,4)& (3,2,6)\\
	 (3,5,0)& (4,0,3)& (4,1,6)& (4,5,2)& (5,6,1)& (6,0,4)& (6,2,3)\end{array}$

	Lexicographic least 7-good sequencing 0123456

	Number of 7-good sequencings found: 197
\end{description}\begin{description}
\item[\D{7}{2}{.57}:]	$\begin{array}{@{}*{7}{l@{\hspace*{6pt}}}@{}}
 (0,6,5)& (1,0,2)& (1,5,3)& (2,0,1)& (2,5,4)& (2,6,3)& (3,1,4)\\
	 (3,5,0)& (3,6,2)& (4,0,3)& (4,5,2)& (4,6,1)& (5,1,6)& (6,0,4)\end{array}$

	Lexicographic least 7-good sequencing 0123456

	Number of 7-good sequencings found: 219
\end{description}\begin{description}
\item[\D{7}{2}{.58}:]	$\begin{array}{@{}*{7}{l@{\hspace*{6pt}}}@{}}
 (0,6,5)& (1,0,2)& (1,5,3)& (2,0,1)& (2,5,4)& (2,6,3)& (3,1,4)\\
	 (3,5,0)& (3,6,2)& (4,0,3)& (4,1,6)& (4,5,2)& (5,6,1)& (6,0,4)\end{array}$

	Lexicographic least 7-good sequencing 0123456

	Number of 7-good sequencings found: 195
\end{description}\begin{description}
\item[\D{7}{2}{.59}:]	$\begin{array}{@{}*{7}{l@{\hspace*{6pt}}}@{}}
 (0,6,5)& (1,0,2)& (1,5,3)& (2,0,1)& (2,4,5)& (3,1,4)& (3,2,6)\\
	 (3,5,0)& (4,0,3)& (4,6,1)& (5,1,6)& (5,4,2)& (6,0,4)& (6,2,3)\end{array}$

	Lexicographic least 7-good sequencing 0123546

	Number of 7-good sequencings found: 128
\end{description}\begin{description}
\item[\D{7}{2}{.60}:]	$\begin{array}{@{}*{7}{l@{\hspace*{6pt}}}@{}}
 (0,6,5)& (1,0,2)& (1,5,3)& (2,0,1)& (2,4,5)& (3,1,4)& (3,2,6)\\
	 (3,5,0)& (4,0,3)& (4,1,6)& (5,4,2)& (5,6,1)& (6,0,4)& (6,2,3)\end{array}$

	Lexicographic least 7-good sequencing 0123546

	Number of 7-good sequencings found: 183
\end{description}\begin{description}
\item[\D{7}{2}{.61}:]	$\begin{array}{@{}*{7}{l@{\hspace*{6pt}}}@{}}
 (0,6,5)& (1,0,2)& (1,5,3)& (2,0,1)& (2,4,5)& (2,6,3)& (3,1,4)\\
	 (3,5,0)& (3,6,2)& (4,0,3)& (4,6,1)& (5,1,6)& (5,4,2)& (6,0,4)\end{array}$

	Lexicographic least 7-good sequencing 0123546

	Number of 7-good sequencings found: 160
\end{description}\begin{description}
\item[\D{7}{2}{.62}:]	$\begin{array}{@{}*{7}{l@{\hspace*{6pt}}}@{}}
 (0,6,5)& (1,0,2)& (1,5,3)& (2,0,1)& (2,4,5)& (2,6,3)& (3,1,4)\\
	 (3,5,0)& (3,6,2)& (4,0,3)& (4,1,6)& (5,4,2)& (5,6,1)& (6,0,4)\end{array}$

	Lexicographic least 7-good sequencing 0123546

	Number of 7-good sequencings found: 169
\end{description}\begin{description}
\item[\D{7}{2}{.63}:]	$\begin{array}{@{}*{7}{l@{\hspace*{6pt}}}@{}}
 (0,6,5)& (1,0,2)& (1,5,3)& (2,0,1)& (2,3,6)& (3,1,4)& (3,5,0)\\
	 (4,0,3)& (4,2,5)& (4,6,1)& (5,1,6)& (5,2,4)& (6,0,4)& (6,3,2)\end{array}$

	Lexicographic least 7-good sequencing 0132456

	Number of 7-good sequencings found: 118
\end{description}\begin{description}
\item[\D{7}{2}{.64}:]	$\begin{array}{@{}*{7}{l@{\hspace*{6pt}}}@{}}
 (0,6,5)& (1,0,2)& (1,5,3)& (2,0,1)& (2,3,6)& (3,1,4)& (3,5,0)\\
	 (4,0,3)& (4,1,6)& (4,2,5)& (5,2,4)& (5,6,1)& (6,0,4)& (6,3,2)\end{array}$

	Lexicographic least 7-good sequencing 0132456

	Number of 7-good sequencings found: 112
\end{description}\begin{description}
\item[\D{7}{2}{.65}:]	$\begin{array}{@{}*{7}{l@{\hspace*{6pt}}}@{}}
 (0,6,5)& (1,0,2)& (1,5,3)& (2,0,1)& (2,3,6)& (2,5,4)& (3,1,4)\\
	 (3,5,0)& (4,0,3)& (4,5,2)& (4,6,1)& (5,1,6)& (6,0,4)& (6,3,2)\end{array}$

	Lexicographic least 7-good sequencing 0132456

	Number of 7-good sequencings found: 178
\end{description}\begin{description}
\item[\D{7}{2}{.66}:]	$\begin{array}{@{}*{7}{l@{\hspace*{6pt}}}@{}}
 (0,6,5)& (1,0,2)& (1,5,3)& (2,0,1)& (2,3,6)& (2,5,4)& (3,1,4)\\
	 (3,5,0)& (4,0,3)& (4,1,6)& (4,5,2)& (5,6,1)& (6,0,4)& (6,3,2)\end{array}$

	Lexicographic least 7-good sequencing 0132456

	Number of 7-good sequencings found: 166
\end{description}\begin{description}
\item[\D{7}{2}{.67}:]	$\begin{array}{@{}*{7}{l@{\hspace*{6pt}}}@{}}
 (0,6,5)& (1,0,2)& (1,5,3)& (2,0,1)& (2,3,6)& (2,4,5)& (3,1,4)\\
	 (3,5,0)& (4,0,3)& (4,6,1)& (5,1,6)& (5,4,2)& (6,0,4)& (6,3,2)\end{array}$

	Lexicographic least 7-good sequencing 0132546

	Number of 7-good sequencings found: 156
\end{description}\begin{description}
\item[\D{7}{2}{.68}:]	$\begin{array}{@{}*{7}{l@{\hspace*{6pt}}}@{}}
 (0,6,5)& (1,0,2)& (1,5,3)& (2,0,1)& (2,3,6)& (2,4,5)& (3,1,4)\\
	 (3,5,0)& (4,0,3)& (4,1,6)& (5,4,2)& (5,6,1)& (6,0,4)& (6,3,2)\end{array}$

	Lexicographic least 7-good sequencing 0132546

	Number of 7-good sequencings found: 138
\end{description}\begin{description}
\item[\D{7}{2}{.69}:]	$\begin{array}{@{}*{7}{l@{\hspace*{6pt}}}@{}}
 (0,6,5)& (1,0,2)& (1,4,6)& (2,0,1)& (3,1,5)& (3,2,6)& (3,4,0)\\
	 (4,1,3)& (4,2,5)& (5,0,3)& (5,2,4)& (5,6,1)& (6,0,4)& (6,2,3)\end{array}$

	Lexicographic least 7-good sequencing 0123564

	Number of 7-good sequencings found: 155
\end{description}\begin{description}
\item[\D{7}{2}{.70}:]	$\begin{array}{@{}*{7}{l@{\hspace*{6pt}}}@{}}
 (0,6,5)& (1,0,2)& (1,4,6)& (2,0,1)& (3,0,4)& (3,1,5)& (3,2,6)\\
	 (4,1,3)& (4,2,5)& (5,0,3)& (5,2,4)& (5,6,1)& (6,2,3)& (6,4,0)\end{array}$

	Lexicographic least 7-good sequencing 0123564

	Number of 7-good sequencings found: 152
\end{description}\begin{description}
\item[\D{7}{2}{.71}:]	$\begin{array}{@{}*{7}{l@{\hspace*{6pt}}}@{}}
 (0,6,5)& (1,0,2)& (1,4,6)& (2,0,1)& (2,6,3)& (3,1,5)& (3,4,0)\\
	 (3,6,2)& (4,1,3)& (4,2,5)& (5,0,3)& (5,2,4)& (5,6,1)& (6,0,4)\end{array}$

	Lexicographic least 7-good sequencing 0123564

	Number of 7-good sequencings found: 175
\end{description}\begin{description}
\item[\D{7}{2}{.72}:]	$\begin{array}{@{}*{7}{l@{\hspace*{6pt}}}@{}}
 (0,6,5)& (1,0,2)& (1,4,6)& (2,0,1)& (2,6,3)& (3,0,4)& (3,1,5)\\
	 (3,6,2)& (4,1,3)& (4,2,5)& (5,0,3)& (5,2,4)& (5,6,1)& (6,4,0)\end{array}$

	Lexicographic least 7-good sequencing 0123564

	Number of 7-good sequencings found: 197
\end{description}\begin{description}
\item[\D{7}{2}{.73}:]	$\begin{array}{@{}*{7}{l@{\hspace*{6pt}}}@{}}
 (0,6,5)& (1,0,2)& (1,4,6)& (2,0,1)& (2,5,4)& (3,1,5)& (3,2,6)\\
	 (3,4,0)& (4,1,3)& (4,5,2)& (5,0,3)& (5,6,1)& (6,0,4)& (6,2,3)\end{array}$

	Lexicographic least 7-good sequencing 0135624

	Number of 7-good sequencings found: 113
\end{description}\begin{description}
\item[\D{7}{2}{.74}:]	$\begin{array}{@{}*{7}{l@{\hspace*{6pt}}}@{}}
 (0,6,5)& (1,0,2)& (1,4,6)& (2,0,1)& (2,5,4)& (3,0,4)& (3,1,5)\\
	 (3,2,6)& (4,1,3)& (4,5,2)& (5,0,3)& (5,6,1)& (6,2,3)& (6,4,0)\end{array}$

	Lexicographic least 7-good sequencing 0135624

	Number of 7-good sequencings found: 112
\end{description}\begin{description}
\item[\D{7}{2}{.75}:]	$\begin{array}{@{}*{7}{l@{\hspace*{6pt}}}@{}}
 (0,6,5)& (1,0,2)& (1,4,6)& (2,0,1)& (2,5,4)& (2,6,3)& (3,1,5)\\
	 (3,4,0)& (3,6,2)& (4,1,3)& (4,5,2)& (5,0,3)& (5,6,1)& (6,0,4)\end{array}$

	Lexicographic least 7-good sequencing 0135264

	Number of 7-good sequencings found: 164
\end{description}\begin{description}
\item[\D{7}{2}{.76}:]	$\begin{array}{@{}*{7}{l@{\hspace*{6pt}}}@{}}
 (0,6,5)& (1,0,2)& (1,4,6)& (2,0,1)& (2,5,4)& (2,6,3)& (3,0,4)\\
	 (3,1,5)& (3,6,2)& (4,1,3)& (4,5,2)& (5,0,3)& (5,6,1)& (6,4,0)\end{array}$

	Lexicographic least 7-good sequencing 0135264

	Number of 7-good sequencings found: 166
\end{description}\begin{description}
\item[\D{7}{2}{.77}:]	$\begin{array}{@{}*{7}{l@{\hspace*{6pt}}}@{}}
 (0,6,5)& (1,0,2)& (1,4,6)& (2,0,1)& (2,4,5)& (3,1,5)& (3,2,6)\\
	 (3,4,0)& (4,1,3)& (5,0,3)& (5,4,2)& (5,6,1)& (6,0,4)& (6,2,3)\end{array}$

	Lexicographic least 7-good sequencing 0123564

	Number of 7-good sequencings found: 150
\end{description}\begin{description}
\item[\D{7}{2}{.78}:]	$\begin{array}{@{}*{7}{l@{\hspace*{6pt}}}@{}}
 (0,6,5)& (1,0,2)& (1,4,6)& (2,0,1)& (2,4,5)& (3,0,4)& (3,1,5)\\
	 (3,2,6)& (4,1,3)& (5,0,3)& (5,4,2)& (5,6,1)& (6,2,3)& (6,4,0)\end{array}$

	Lexicographic least 7-good sequencing 0123564

	Number of 7-good sequencings found: 148
\end{description}\begin{description}
\item[\D{7}{2}{.79}:]	$\begin{array}{@{}*{7}{l@{\hspace*{6pt}}}@{}}
 (0,6,5)& (1,0,2)& (1,4,6)& (2,0,1)& (2,4,5)& (2,6,3)& (3,1,5)\\
	 (3,4,0)& (3,6,2)& (4,1,3)& (5,0,3)& (5,4,2)& (5,6,1)& (6,0,4)\end{array}$

	Lexicographic least 7-good sequencing 0123564

	Number of 7-good sequencings found: 201
\end{description}\begin{description}
\item[\D{7}{2}{.80}:]	$\begin{array}{@{}*{7}{l@{\hspace*{6pt}}}@{}}
 (0,6,5)& (1,0,2)& (1,4,6)& (2,0,1)& (2,4,5)& (2,6,3)& (3,0,4)\\
	 (3,1,5)& (3,6,2)& (4,1,3)& (5,0,3)& (5,4,2)& (5,6,1)& (6,4,0)\end{array}$

	Lexicographic least 7-good sequencing 0123564

	Number of 7-good sequencings found: 195
\end{description}\begin{description}
\item[\D{7}{2}{.81}:]	$\begin{array}{@{}*{7}{l@{\hspace*{6pt}}}@{}}
 (0,6,5)& (1,0,2)& (1,4,6)& (2,0,1)& (2,3,6)& (3,1,5)& (3,4,0)\\
	 (4,1,3)& (4,2,5)& (5,0,3)& (5,2,4)& (5,6,1)& (6,0,4)& (6,3,2)\end{array}$

	Lexicographic least 7-good sequencing 0125634

	Number of 7-good sequencings found: 160
\end{description}\begin{description}
\item[\D{7}{2}{.82}:]	$\begin{array}{@{}*{7}{l@{\hspace*{6pt}}}@{}}
 (0,6,5)& (1,0,2)& (1,4,6)& (2,0,1)& (2,3,6)& (3,0,4)& (3,1,5)\\
	 (4,1,3)& (4,2,5)& (5,0,3)& (5,2,4)& (5,6,1)& (6,3,2)& (6,4,0)\end{array}$

	Lexicographic least 7-good sequencing 0125634

	Number of 7-good sequencings found: 183
\end{description}\begin{description}
\item[\D{7}{2}{.83}:]	$\begin{array}{@{}*{7}{l@{\hspace*{6pt}}}@{}}
 (0,6,5)& (1,0,2)& (1,4,6)& (2,0,1)& (2,3,6)& (2,5,4)& (3,1,5)\\
	 (3,4,0)& (4,1,3)& (4,5,2)& (5,0,3)& (5,6,1)& (6,0,4)& (6,3,2)\end{array}$

	Lexicographic least 7-good sequencing 0135264

	Number of 7-good sequencings found: 163
\end{description}\begin{description}
\item[\D{7}{2}{.84}:]	$\begin{array}{@{}*{7}{l@{\hspace*{6pt}}}@{}}
 (0,6,5)& (1,0,2)& (1,4,6)& (2,0,1)& (2,3,6)& (2,5,4)& (3,0,4)\\
	 (3,1,5)& (4,1,3)& (4,5,2)& (5,0,3)& (5,6,1)& (6,3,2)& (6,4,0)\end{array}$

	Lexicographic least 7-good sequencing 0135264

	Number of 7-good sequencings found: 138
\end{description}\begin{description}
\item[\D{7}{2}{.85}:]	$\begin{array}{@{}*{7}{l@{\hspace*{6pt}}}@{}}
 (0,6,5)& (1,0,2)& (1,4,6)& (2,0,1)& (2,3,6)& (2,4,5)& (3,1,5)\\
	 (3,4,0)& (4,1,3)& (5,0,3)& (5,4,2)& (5,6,1)& (6,0,4)& (6,3,2)\end{array}$

	Lexicographic least 7-good sequencing 0125634

	Number of 7-good sequencings found: 163
\end{description}\begin{description}
\item[\D{7}{2}{.86}:]	$\begin{array}{@{}*{7}{l@{\hspace*{6pt}}}@{}}
 (0,6,5)& (1,0,2)& (1,4,6)& (2,0,1)& (2,3,6)& (2,4,5)& (3,0,4)\\
	 (3,1,5)& (4,1,3)& (5,0,3)& (5,4,2)& (5,6,1)& (6,3,2)& (6,4,0)\end{array}$

	Lexicographic least 7-good sequencing 0125634

	Number of 7-good sequencings found: 169
\end{description}\begin{description}
\item[\D{7}{2}{.87}:]	$\begin{array}{@{}*{7}{l@{\hspace*{6pt}}}@{}}
 (0,6,5)& (1,0,2)& (1,4,6)& (1,5,3)& (2,0,1)& (3,2,6)& (3,4,1)\\
	 (3,5,0)& (4,0,3)& (4,2,5)& (5,2,4)& (5,6,1)& (6,0,4)& (6,2,3)\end{array}$

	Lexicographic least 7-good sequencing 0123564

	Number of 7-good sequencings found: 177
\end{description}\begin{description}
\item[\D{7}{2}{.88}:]	$\begin{array}{@{}*{7}{l@{\hspace*{6pt}}}@{}}
 (0,6,5)& (1,0,2)& (1,4,6)& (1,5,3)& (2,0,1)& (2,6,3)& (3,4,1)\\
	 (3,5,0)& (3,6,2)& (4,0,3)& (4,2,5)& (5,2,4)& (5,6,1)& (6,0,4)\end{array}$

	Lexicographic least 7-good sequencing 0123564

	Number of 7-good sequencings found: 180
\end{description}\begin{description}
\item[\D{7}{2}{.89}:]	$\begin{array}{@{}*{7}{l@{\hspace*{6pt}}}@{}}
 (0,6,5)& (1,0,2)& (1,4,6)& (1,5,3)& (2,0,1)& (2,5,4)& (3,2,6)\\
	 (3,4,1)& (3,5,0)& (4,0,3)& (4,5,2)& (5,6,1)& (6,0,4)& (6,2,3)\end{array}$

	Lexicographic least 7-good sequencing 0135624

	Number of 7-good sequencings found: 210
\end{description}\begin{description}
\item[\D{7}{2}{.90}:]	$\begin{array}{@{}*{7}{l@{\hspace*{6pt}}}@{}}
 (0,6,5)& (1,0,2)& (1,4,6)& (1,5,3)& (2,0,1)& (2,5,4)& (2,6,3)\\
	 (3,4,1)& (3,5,0)& (3,6,2)& (4,0,3)& (4,5,2)& (5,6,1)& (6,0,4)\end{array}$

	Lexicographic least 7-good sequencing 0135264

	Number of 7-good sequencings found: 214
\end{description}\begin{description}
\item[\D{7}{2}{.91}:]	$\begin{array}{@{}*{7}{l@{\hspace*{6pt}}}@{}}
 (0,6,5)& (1,0,2)& (1,4,6)& (1,5,3)& (2,0,1)& (2,4,5)& (3,2,6)\\
	 (3,4,1)& (3,5,0)& (4,0,3)& (5,4,2)& (5,6,1)& (6,0,4)& (6,2,3)\end{array}$

	Lexicographic least 7-good sequencing 0123564

	Number of 7-good sequencings found: 191
\end{description}\begin{description}
\item[\D{7}{2}{.92}:]	$\begin{array}{@{}*{7}{l@{\hspace*{6pt}}}@{}}
 (0,6,5)& (1,0,2)& (1,4,6)& (1,5,3)& (2,0,1)& (2,3,6)& (3,4,1)\\
	 (3,5,0)& (4,0,3)& (4,2,5)& (5,2,4)& (5,6,1)& (6,0,4)& (6,3,2)\end{array}$

	Lexicographic least 7-good sequencing 0132564

	Number of 7-good sequencings found: 157
\end{description}\begin{description}
\item[\D{7}{2}{.93}:]	$\begin{array}{@{}*{7}{l@{\hspace*{6pt}}}@{}}
 (0,6,5)& (1,0,2)& (1,4,6)& (1,5,3)& (2,0,1)& (2,3,6)& (2,5,4)\\
	 (3,4,1)& (3,5,0)& (4,0,3)& (4,5,2)& (5,6,1)& (6,0,4)& (6,3,2)\end{array}$

	Lexicographic least 7-good sequencing 0135264

	Number of 7-good sequencings found: 196
\end{description}\begin{description}
\item[\D{7}{2}{.94}:]	$\begin{array}{@{}*{7}{l@{\hspace*{6pt}}}@{}}
 (0,6,5)& (1,0,2)& (1,4,6)& (1,5,3)& (2,0,1)& (2,3,6)& (2,4,5)\\
	 (3,4,1)& (3,5,0)& (4,0,3)& (5,4,2)& (5,6,1)& (6,0,4)& (6,3,2)\end{array}$

	Lexicographic least 7-good sequencing 0132564

	Number of 7-good sequencings found: 181
\end{description}\begin{description}
\item[\D{7}{2}{.95}:]	$\begin{array}{@{}*{7}{l@{\hspace*{6pt}}}@{}}
 (0,6,5)& (1,0,2)& (1,4,3)& (2,0,1)& (3,1,5)& (3,2,6)& (3,4,0)\\
	 (4,2,5)& (4,6,1)& (5,0,3)& (5,1,6)& (5,2,4)& (6,0,4)& (6,2,3)\end{array}$

	Lexicographic least 7-good sequencing 0123456

	Number of 7-good sequencings found: 155
\end{description}\begin{description}
\item[\D{7}{2}{.96}:]	$\begin{array}{@{}*{7}{l@{\hspace*{6pt}}}@{}}
 (0,6,5)& (1,0,2)& (1,4,3)& (2,0,1)& (3,1,5)& (3,2,6)& (3,4,0)\\
	 (4,1,6)& (4,2,5)& (5,0,3)& (5,2,4)& (5,6,1)& (6,0,4)& (6,2,3)\end{array}$

	Lexicographic least 7-good sequencing 0123456

	Number of 7-good sequencings found: 140
\end{description}\begin{description}
\item[\D{7}{2}{.97}:]	$\begin{array}{@{}*{7}{l@{\hspace*{6pt}}}@{}}
 (0,6,5)& (1,0,2)& (1,4,3)& (2,0,1)& (3,0,4)& (3,1,5)& (3,2,6)\\
	 (4,2,5)& (4,6,1)& (5,0,3)& (5,1,6)& (5,2,4)& (6,2,3)& (6,4,0)\end{array}$

	Lexicographic least 7-good sequencing 0123456

	Number of 7-good sequencings found: 148
\end{description}\begin{description}
\item[\D{7}{2}{.98}:]	$\begin{array}{@{}*{7}{l@{\hspace*{6pt}}}@{}}
 (0,6,5)& (1,0,2)& (1,4,3)& (2,0,1)& (3,0,4)& (3,1,5)& (3,2,6)\\
	 (4,2,5)& (4,6,0)& (5,0,3)& (5,1,6)& (5,2,4)& (6,2,3)& (6,4,1)\end{array}$

	Lexicographic least 7-good sequencing 0123456

	Number of 7-good sequencings found: 206
\end{description}\begin{description}
\item[\D{7}{2}{.99}:]	$\begin{array}{@{}*{7}{l@{\hspace*{6pt}}}@{}}
 (0,6,5)& (1,0,2)& (1,4,3)& (2,0,1)& (3,0,4)& (3,1,5)& (3,2,6)\\
	 (4,1,6)& (4,2,5)& (5,0,3)& (5,2,4)& (5,6,1)& (6,2,3)& (6,4,0)\end{array}$

	Lexicographic least 7-good sequencing 0123456

	Number of 7-good sequencings found: 153
\end{description}\begin{description}
\item[\D{7}{2}{.100}:]	$\begin{array}{@{}*{7}{l@{\hspace*{6pt}}}@{}}
 (0,6,5)& (1,0,2)& (1,4,3)& (2,0,1)& (2,6,3)& (3,1,5)& (3,4,0)\\
	 (3,6,2)& (4,2,5)& (4,6,1)& (5,0,3)& (5,1,6)& (5,2,4)& (6,0,4)\end{array}$

	Lexicographic least 7-good sequencing 0123456

	Number of 7-good sequencings found: 150
\end{description}\begin{description}
\item[\D{7}{2}{.101}:]	$\begin{array}{@{}*{7}{l@{\hspace*{6pt}}}@{}}
 (0,6,5)& (1,0,2)& (1,4,3)& (2,0,1)& (2,6,3)& (3,1,5)& (3,4,0)\\
	 (3,6,2)& (4,1,6)& (4,2,5)& (5,0,3)& (5,2,4)& (5,6,1)& (6,0,4)\end{array}$

	Lexicographic least 7-good sequencing 0123456

	Number of 7-good sequencings found: 146
\end{description}\begin{description}
\item[\D{7}{2}{.102}:]	$\begin{array}{@{}*{7}{l@{\hspace*{6pt}}}@{}}
 (0,6,5)& (1,0,2)& (1,4,3)& (2,0,1)& (2,6,3)& (3,0,4)& (3,1,5)\\
	 (3,6,2)& (4,2,5)& (4,6,1)& (5,0,3)& (5,1,6)& (5,2,4)& (6,4,0)\end{array}$

	Lexicographic least 7-good sequencing 0123456

	Number of 7-good sequencings found: 168
\end{description}\begin{description}
\item[\D{7}{2}{.103}:]	$\begin{array}{@{}*{7}{l@{\hspace*{6pt}}}@{}}
 (0,6,5)& (1,0,2)& (1,4,3)& (2,0,1)& (2,6,3)& (3,0,4)& (3,1,5)\\
	 (3,6,2)& (4,1,6)& (4,2,5)& (5,0,3)& (5,2,4)& (5,6,1)& (6,4,0)\end{array}$

	Lexicographic least 7-good sequencing 0123456

	Number of 7-good sequencings found: 177
\end{description}\begin{description}
\item[\D{7}{2}{.104}:]	$\begin{array}{@{}*{7}{l@{\hspace*{6pt}}}@{}}
 (0,6,5)& (1,0,2)& (1,4,3)& (2,0,1)& (2,5,4)& (3,1,5)& (3,2,6)\\
	 (3,4,0)& (4,5,2)& (4,6,1)& (5,0,3)& (5,1,6)& (6,0,4)& (6,2,3)\end{array}$

	Lexicographic least 7-good sequencing 0123456

	Number of 7-good sequencings found: 160
\end{description}\begin{description}
\item[\D{7}{2}{.105}:]	$\begin{array}{@{}*{7}{l@{\hspace*{6pt}}}@{}}
 (0,6,5)& (1,0,2)& (1,4,3)& (2,0,1)& (2,5,4)& (3,1,5)& (3,2,6)\\
	 (3,4,0)& (4,1,6)& (4,5,2)& (5,0,3)& (5,6,1)& (6,0,4)& (6,2,3)\end{array}$

	Lexicographic least 7-good sequencing 0123456

	Number of 7-good sequencings found: 98
\end{description}\begin{description}
\item[\D{7}{2}{.106}:]	$\begin{array}{@{}*{7}{l@{\hspace*{6pt}}}@{}}
 (0,6,5)& (1,0,2)& (1,4,3)& (2,0,1)& (2,5,4)& (3,0,4)& (3,1,5)\\
	 (3,2,6)& (4,5,2)& (4,6,1)& (5,0,3)& (5,1,6)& (6,2,3)& (6,4,0)\end{array}$

	Lexicographic least 7-good sequencing 0123456

	Number of 7-good sequencings found: 158
\end{description}\begin{description}
\item[\D{7}{2}{.107}:]	$\begin{array}{@{}*{7}{l@{\hspace*{6pt}}}@{}}
 (0,6,5)& (1,0,2)& (1,4,3)& (2,0,1)& (2,5,4)& (3,0,4)& (3,1,5)\\
	 (3,2,6)& (4,5,2)& (4,6,0)& (5,0,3)& (5,1,6)& (6,2,3)& (6,4,1)\end{array}$

	Lexicographic least 7-good sequencing 0123456

	Number of 7-good sequencings found: 184
\end{description}\begin{description}
\item[\D{7}{2}{.108}:]	$\begin{array}{@{}*{7}{l@{\hspace*{6pt}}}@{}}
 (0,6,5)& (1,0,2)& (1,4,3)& (2,0,1)& (2,5,4)& (3,0,4)& (3,1,5)\\
	 (3,2,6)& (4,1,6)& (4,5,2)& (5,0,3)& (5,6,1)& (6,2,3)& (6,4,0)\end{array}$

	Lexicographic least 7-good sequencing 0123456

	Number of 7-good sequencings found: 118
\end{description}\begin{description}
\item[\D{7}{2}{.109}:]	$\begin{array}{@{}*{7}{l@{\hspace*{6pt}}}@{}}
 (0,6,5)& (1,0,2)& (1,4,3)& (2,0,1)& (2,5,4)& (2,6,3)& (3,1,5)\\
	 (3,4,0)& (3,6,2)& (4,5,2)& (4,6,1)& (5,0,3)& (5,1,6)& (6,0,4)\end{array}$

	Lexicographic least 7-good sequencing 0123456

	Number of 7-good sequencings found: 163
\end{description}\begin{description}
\item[\D{7}{2}{.110}:]	$\begin{array}{@{}*{7}{l@{\hspace*{6pt}}}@{}}
 (0,6,5)& (1,0,2)& (1,4,3)& (2,0,1)& (2,5,4)& (2,6,3)& (3,1,5)\\
	 (3,4,0)& (3,6,2)& (4,1,6)& (4,5,2)& (5,0,3)& (5,6,1)& (6,0,4)\end{array}$

	Lexicographic least 7-good sequencing 0123456

	Number of 7-good sequencings found: 152
\end{description}\begin{description}
\item[\D{7}{2}{.111}:]	$\begin{array}{@{}*{7}{l@{\hspace*{6pt}}}@{}}
 (0,6,5)& (1,0,2)& (1,4,3)& (2,0,1)& (2,5,4)& (2,6,3)& (3,0,4)\\
	 (3,1,5)& (3,6,2)& (4,5,2)& (4,6,1)& (5,0,3)& (5,1,6)& (6,4,0)\end{array}$

	Lexicographic least 7-good sequencing 0123456

	Number of 7-good sequencings found: 179
\end{description}\begin{description}
\item[\D{7}{2}{.112}:]	$\begin{array}{@{}*{7}{l@{\hspace*{6pt}}}@{}}
 (0,6,5)& (1,0,2)& (1,4,3)& (2,0,1)& (2,5,4)& (2,6,3)& (3,0,4)\\
	 (3,1,5)& (3,6,2)& (4,1,6)& (4,5,2)& (5,0,3)& (5,6,1)& (6,4,0)\end{array}$

	Lexicographic least 7-good sequencing 0123456

	Number of 7-good sequencings found: 178
\end{description}\begin{description}
\item[\D{7}{2}{.113}:]	$\begin{array}{@{}*{7}{l@{\hspace*{6pt}}}@{}}
 (0,6,5)& (1,0,2)& (1,4,3)& (2,0,1)& (2,4,5)& (3,1,5)& (3,2,6)\\
	 (3,4,0)& (4,6,1)& (5,0,3)& (5,1,6)& (5,4,2)& (6,0,4)& (6,2,3)\end{array}$

	Lexicographic least 7-good sequencing 0123546

	Number of 7-good sequencings found: 175
\end{description}\begin{description}
\item[\D{7}{2}{.114}:]	$\begin{array}{@{}*{7}{l@{\hspace*{6pt}}}@{}}
 (0,6,5)& (1,0,2)& (1,4,3)& (2,0,1)& (2,4,5)& (3,1,5)& (3,2,6)\\
	 (3,4,0)& (4,1,6)& (5,0,3)& (5,4,2)& (5,6,1)& (6,0,4)& (6,2,3)\end{array}$

	Lexicographic least 7-good sequencing 0123546

	Number of 7-good sequencings found: 146
\end{description}\begin{description}
\item[\D{7}{2}{.115}:]	$\begin{array}{@{}*{7}{l@{\hspace*{6pt}}}@{}}
 (0,6,5)& (1,0,2)& (1,4,3)& (2,0,1)& (2,4,5)& (3,0,4)& (3,1,5)\\
	 (3,2,6)& (4,6,1)& (5,0,3)& (5,1,6)& (5,4,2)& (6,2,3)& (6,4,0)\end{array}$

	Lexicographic least 7-good sequencing 0123546

	Number of 7-good sequencings found: 166
\end{description}\begin{description}
\item[\D{7}{2}{.116}:]	$\begin{array}{@{}*{7}{l@{\hspace*{6pt}}}@{}}
 (0,6,5)& (1,0,2)& (1,4,3)& (2,0,1)& (2,4,5)& (3,0,4)& (3,1,5)\\
	 (3,2,6)& (4,1,6)& (5,0,3)& (5,4,2)& (5,6,1)& (6,2,3)& (6,4,0)\end{array}$

	Lexicographic least 7-good sequencing 0123546

	Number of 7-good sequencings found: 155
\end{description}\begin{description}
\item[\D{7}{2}{.117}:]	$\begin{array}{@{}*{7}{l@{\hspace*{6pt}}}@{}}
 (0,6,5)& (1,0,2)& (1,4,3)& (2,0,1)& (2,4,5)& (2,6,3)& (3,1,5)\\
	 (3,4,0)& (3,6,2)& (4,6,1)& (5,0,3)& (5,1,6)& (5,4,2)& (6,0,4)\end{array}$

	Lexicographic least 7-good sequencing 0123546

	Number of 7-good sequencings found: 201
\end{description}\begin{description}
\item[\D{7}{2}{.118}:]	$\begin{array}{@{}*{7}{l@{\hspace*{6pt}}}@{}}
 (0,6,5)& (1,0,2)& (1,4,3)& (2,0,1)& (2,4,5)& (2,6,3)& (3,1,5)\\
	 (3,4,0)& (3,6,2)& (4,1,6)& (5,0,3)& (5,4,2)& (5,6,1)& (6,0,4)\end{array}$

	Lexicographic least 7-good sequencing 0123546

	Number of 7-good sequencings found: 228
\end{description}\begin{description}
\item[\D{7}{2}{.119}:]	$\begin{array}{@{}*{7}{l@{\hspace*{6pt}}}@{}}
 (0,6,5)& (1,0,2)& (1,4,3)& (2,0,1)& (2,4,5)& (2,6,3)& (3,0,4)\\
	 (3,1,5)& (3,6,2)& (4,6,1)& (5,0,3)& (5,1,6)& (5,4,2)& (6,4,0)\end{array}$

	Lexicographic least 7-good sequencing 0123546

	Number of 7-good sequencings found: 195
\end{description}\begin{description}
\item[\D{7}{2}{.120}:]	$\begin{array}{@{}*{7}{l@{\hspace*{6pt}}}@{}}
 (0,6,5)& (1,0,2)& (1,4,3)& (2,0,1)& (2,4,5)& (2,6,3)& (3,0,4)\\
	 (3,1,5)& (3,6,2)& (4,1,6)& (5,0,3)& (5,4,2)& (5,6,1)& (6,4,0)\end{array}$

	Lexicographic least 7-good sequencing 0123546

	Number of 7-good sequencings found: 219
\end{description}\begin{description}
\item[\D{7}{2}{.121}:]	$\begin{array}{@{}*{7}{l@{\hspace*{6pt}}}@{}}
 (0,6,5)& (1,0,2)& (1,4,3)& (2,0,1)& (2,3,6)& (3,1,5)& (3,4,0)\\
	 (4,2,5)& (4,6,1)& (5,0,3)& (5,1,6)& (5,2,4)& (6,0,4)& (6,3,2)\end{array}$

	Lexicographic least 7-good sequencing 0125634

	Number of 7-good sequencings found: 113
\end{description}\begin{description}
\item[\D{7}{2}{.122}:]	$\begin{array}{@{}*{7}{l@{\hspace*{6pt}}}@{}}
 (0,6,5)& (1,0,2)& (1,4,3)& (2,0,1)& (2,3,6)& (3,1,5)& (3,4,0)\\
	 (4,1,6)& (4,2,5)& (5,0,3)& (5,2,4)& (5,6,1)& (6,0,4)& (6,3,2)\end{array}$

	Lexicographic least 7-good sequencing 0125634

	Number of 7-good sequencings found: 98
\end{description}\begin{description}
\item[\D{7}{2}{.123}:]	$\begin{array}{@{}*{7}{l@{\hspace*{6pt}}}@{}}
 (0,6,5)& (1,0,2)& (1,4,3)& (2,0,1)& (2,3,6)& (3,0,4)& (3,1,5)\\
	 (4,2,5)& (4,6,1)& (5,0,3)& (5,1,6)& (5,2,4)& (6,3,2)& (6,4,0)\end{array}$

	Lexicographic least 7-good sequencing 0125634

	Number of 7-good sequencings found: 122
\end{description}\begin{description}
\item[\D{7}{2}{.124}:]	$\begin{array}{@{}*{7}{l@{\hspace*{6pt}}}@{}}
 (0,6,5)& (1,0,2)& (1,4,3)& (2,0,1)& (2,3,6)& (3,0,4)& (3,1,5)\\
	 (4,2,5)& (4,6,0)& (5,0,3)& (5,1,6)& (5,2,4)& (6,3,2)& (6,4,1)\end{array}$

	Lexicographic least 7-good sequencing 0125634

	Number of 7-good sequencings found: 189
\end{description}\begin{description}
\item[\D{7}{2}{.125}:]	$\begin{array}{@{}*{7}{l@{\hspace*{6pt}}}@{}}
 (0,6,5)& (1,0,2)& (1,4,3)& (2,0,1)& (2,3,6)& (3,0,4)& (3,1,5)\\
	 (4,1,6)& (4,2,5)& (5,0,3)& (5,2,4)& (5,6,1)& (6,3,2)& (6,4,0)\end{array}$

	Lexicographic least 7-good sequencing 0125634

	Number of 7-good sequencings found: 128
\end{description}\begin{description}
\item[\D{7}{2}{.126}:]	$\begin{array}{@{}*{7}{l@{\hspace*{6pt}}}@{}}
 (0,6,5)& (1,0,2)& (1,4,3)& (2,0,1)& (2,3,6)& (2,5,4)& (3,1,5)\\
	 (3,4,0)& (4,5,2)& (4,6,1)& (5,0,3)& (5,1,6)& (6,0,4)& (6,3,2)\end{array}$

	Lexicographic least 7-good sequencing 0132456

	Number of 7-good sequencings found: 163
\end{description}\begin{description}
\item[\D{7}{2}{.127}:]	$\begin{array}{@{}*{7}{l@{\hspace*{6pt}}}@{}}
 (0,6,5)& (1,0,2)& (1,4,3)& (2,0,1)& (2,3,6)& (2,5,4)& (3,1,5)\\
	 (3,4,0)& (4,1,6)& (4,5,2)& (5,0,3)& (5,6,1)& (6,0,4)& (6,3,2)\end{array}$

	Lexicographic least 7-good sequencing 0132456

	Number of 7-good sequencings found: 160
\end{description}\begin{description}
\item[\D{7}{2}{.128}:]	$\begin{array}{@{}*{7}{l@{\hspace*{6pt}}}@{}}
 (0,6,5)& (1,0,2)& (1,4,3)& (2,0,1)& (2,3,6)& (2,5,4)& (3,0,4)\\
	 (3,1,5)& (4,5,2)& (4,6,1)& (5,0,3)& (5,1,6)& (6,3,2)& (6,4,0)\end{array}$

	Lexicographic least 7-good sequencing 0132456

	Number of 7-good sequencings found: 139
\end{description}\begin{description}
\item[\D{7}{2}{.129}:]	$\begin{array}{@{}*{7}{l@{\hspace*{6pt}}}@{}}
 (0,6,5)& (1,0,2)& (1,4,3)& (2,0,1)& (2,3,6)& (2,5,4)& (3,0,4)\\
	 (3,1,5)& (4,5,2)& (4,6,0)& (5,0,3)& (5,1,6)& (6,3,2)& (6,4,1)\end{array}$

	Lexicographic least 7-good sequencing 0132456

	Number of 7-good sequencings found: 149
\end{description}\begin{description}
\item[\D{7}{2}{.130}:]	$\begin{array}{@{}*{7}{l@{\hspace*{6pt}}}@{}}
 (0,6,5)& (1,0,2)& (1,4,3)& (2,0,1)& (2,3,6)& (2,5,4)& (3,0,4)\\
	 (3,1,5)& (4,1,6)& (4,5,2)& (5,0,3)& (5,6,1)& (6,3,2)& (6,4,0)\end{array}$

	Lexicographic least 7-good sequencing 0132456

	Number of 7-good sequencings found: 156
\end{description}\begin{description}
\item[\D{7}{2}{.131}:]	$\begin{array}{@{}*{7}{l@{\hspace*{6pt}}}@{}}
 (0,6,5)& (1,0,2)& (1,4,3)& (2,0,1)& (2,3,6)& (2,4,5)& (3,1,5)\\
	 (3,4,0)& (4,6,1)& (5,0,3)& (5,1,6)& (5,4,2)& (6,0,4)& (6,3,2)\end{array}$

	Lexicographic least 7-good sequencing 0125634

	Number of 7-good sequencings found: 164
\end{description}\begin{description}
\item[\D{7}{2}{.132}:]	$\begin{array}{@{}*{7}{l@{\hspace*{6pt}}}@{}}
 (0,6,5)& (1,0,2)& (1,4,3)& (2,0,1)& (2,3,6)& (2,4,5)& (3,1,5)\\
	 (3,4,0)& (4,1,6)& (5,0,3)& (5,4,2)& (5,6,1)& (6,0,4)& (6,3,2)\end{array}$

	Lexicographic least 7-good sequencing 0125634

	Number of 7-good sequencings found: 152
\end{description}\begin{description}
\item[\D{7}{2}{.133}:]	$\begin{array}{@{}*{7}{l@{\hspace*{6pt}}}@{}}
 (0,6,5)& (1,0,2)& (1,4,3)& (2,0,1)& (2,3,6)& (2,4,5)& (3,0,4)\\
	 (3,1,5)& (4,6,1)& (5,0,3)& (5,1,6)& (5,4,2)& (6,3,2)& (6,4,0)\end{array}%
$

	Lexicographic least 7-good sequencing 0124356

	Number of 7-good sequencings found: 200
\end{description}

\section*{\boldmath The 1,060 nonisomorphic \DTS$(7)$ with underlying \TTS$(7)$ that is \T{7}{3}{} and the 7-good sequencings found.}

\begin{description}
\item[\D{7}{3}{.1}:]	$%
$

	Lexicographic least 7-good sequencing 0614325

	Number of 7-good sequencings found: 123
\end{description}

\section*{\boldmath The 1,016 nonisomorphic \DTS$(7)$ with underlying \TTS$(7)$ that is \T{7}{4}{} and the 7-good sequencings found.}
\begin{description}
\item[\D{7}{4}{.1}:]	$%
$

	Number of 7-good sequencings found: 0

        Lexicographic least 6-good sequencing 0123456

        Number of 6-good sequencings found: 124
\end{description}\begin{description}
\item[\D{7}{4}{.927}:]	$%
$

	Number of 7-good sequencings found: 0

	Lexicographic least 6-good sequencing 0245613

        Number of 6-good sequencings found: 124
\end{description}\begin{description}
\item[\D{7}{4}{.959}:]	$%
$

	Number of 7-good sequencings found: 0

        Lexicographic least 6-good sequencing 0153462

        Number of 6-good sequencings found: 112

\end{description}\begin{description}
\item[\D{7}{4}{.1016}:]	$%
$

	Number of 7-good sequencings found: 0

        Lexicographic least 6-good sequencing 0124356

        Number of 6-good sequencings found: 112

\end{description}


\begin{thebibliography}{X}

\bibitem{KSV}
D.L. Kreher,  D.R. Stinson and S. Veitch,
Block-avoiding point sequencings of directed triple systems,
\textsl{Preprint} (2019).


\bibitem{OP}
P.R.J. \"{O}sterg\r{a}rd and O. Pottone,
Classification of directed and hybrid triple systems,
\textsl{Bayreuth. Math. Schr.} \textbf{74} (2005), 276--291.

\end{thebibliography}
\end{document}